\documentclass[a4paper,12pt]{amsart}

\usepackage{amssymb}
\usepackage{amsmath}
\usepackage{amsthm}

\usepackage{mathrsfs}

\allowdisplaybreaks
\numberwithin{equation}{section}

\theoremstyle{plain}
\newtheorem{thm}{Theorem}[section]
\newtheorem{prop}[thm]{Proposition}
\newtheorem{cor}[thm]{Corollary}
\newtheorem{lem}[thm]{Lemma}

\theoremstyle{definition}
\newtheorem{defn}[thm]{Definition}
\newtheorem{rem}[thm]{Remark}

\newcommand{\be}{\begin{equation}}
\newcommand{\ee}{\end{equation}}
\def\R{{\mathbb R}}
\def\N{{\mathbb N}}
\def\Z{{\mathbb Z}}
\def\C{{\mathbb C}}
\def\S3{{{\mathbb S}^3}}
\def\SU2{{{\rm SU}(2)}}
\def\Rn{{{\mathbb R}^n}}
\def\Tn{{{\mathbb T}^n}}
\def\Zn{{{\mathbb Z}^n}}
\def\Sn{{{\mathbb S}^n}}

\def\irm{{\rm i}}

\def\drm{{\ {\rm d}}}

\def\Lcal{{\mathcal L}}

\def\Rcal{{\mathcal R}}

\def\Hcal{{\mathcal H}}

\def\Ucal{{\mathcal U}}

\def\va{\varphi}
\def\Lap{{\mathscr L}}
\def\Rep{{{\rm Rep}}}

\def\p#1{{\left({#1}\right)}}

\def\jp#1{{\left\langle{#1}\right\rangle}}
\def\n#1{{\left\|{#1}\right\|}}

\begin{document}

\title[Quantization of pseudo-differential operators
on $\SU2$ and $\S3$]
{Global quantization of pseudo-differential operators
on compact Lie groups, $\SU2$ and 3-sphere}

\author[Michael Ruzhansky]{Michael Ruzhansky}
\address{
  Michael Ruzhansky:
  \endgraf
  Department of Mathematics
  \endgraf
  Imperial College London
  \endgraf
  180 Queen's Gate, London SW7 2AZ 
  \endgraf
  United Kingdom
  \endgraf
  {\it E-mail address} {\rm m.ruzhansky@imperial.ac.uk}
  }

\author[Ville Turunen]{Ville Turunen}
\address{
  Ville Turunen:
  \endgraf
   Helsinki University of Technology
  \endgraf
  Institute of Mathematics
  \endgraf
   P.O. Box 1100
  \endgraf
   FI-02015 HUT
  \endgraf
  Finland
  \endgraf
  {\it E-mail address} {\rm ville.turunen@hut.fi}
  }

\thanks{The first
 author was supported in part by the EPSRC
 Leadership Fellowship EP/G007233/1.}
\date{\today}


\begin{abstract}
Global quantization of 
pseudo-differential operators on 
compact Lie groups is introduced relying on the 
representation theory of the group rather than on expressions
in local coordinates. 
Operators on the 3-dimensional sphere
$\S3$ and on group $\SU2$ are analysed in
detail. A new class of globally defined
symbols is introduced giving rise to the usual H\"ormander's
classes of operators $\Psi^m(G)$, $\Psi^m(\S3)$ and $\Psi^m(\SU2)$.
Properties of the new class and symbolic 
calculus are analysed. 
Properties of symbols as well as $L^2$-boundedness
and Sobolev $L^2$--boundedness
of operators in this global
quantization are established on general compact Lie
groups.
\end{abstract}

\maketitle

\tableofcontents

\section{Introduction}

In this paper we investigate a global quantization of 
operators on compact Lie groups.
We develop a non-commutative analogue of the 
Kohn--Nirenberg quantization of pseudo-differential operators
(\cite{KohnNirenberg}).
The introduced matrix-valued full symbols turn out to have a
number of unexpected properties.
Among other things, the introduced approach provides 
a characterization of the H\"ormander's class 
of pseudo-differential operators on compact Lie
groups using a global
quantization of operators
relying on the representation theory rather than
on the usual expressions in local coordinate charts.
The cases of the 3-dimensional sphere
$\S3$ and Lie group $\SU2$ are analysed in detail
and we show that  pseudo-differential operators
from H\"ormander's classes $\Psi^m$ on these spaces
have matrix-valued
symbols with a remarkable rapid off-diagonal decay property.

There have been many works aiming at the understanding of
pseudo-differential operators on Lie groups, see e.g.
work on left-invariant operators 
\cite{St72, Me83, GK97},
convolution calculus on nilpotent Lie groups 
\cite{MeSh87}, $L^2$-boundedness of
convolution operators related to
the Howe's conjecture \cite{How84,Gl04}, and many others. 
In particular, Theorem \ref{THM:su2-L2} allows $x$-dependence 
and also removes the decay condition on the symbol in the
setting of general compact Lie groups
(a possibility of the
relaxation of decay conditions for derivatives of
symbols with respect to the dual variable for the
$L^2$-boundedness was conjectured in
\cite{How84}).

The present research is inspired by M. Taylor's work
\cite{Taylor84}, who used the exponential mapping to
rely on pseudo-differential operators on the Lie algebra
which can be viewed as the Euclidean space with the
corresponding standard theory of pseudo-differential operators.
However,
the approach developed in this paper is different 
from that of \cite{Taylor84,Taylor86}
since 
it relies on the group structure directly and thus we
do not need to work in neighbourhoods of the neutral
element and can approach global symbol classes directly. 

As usual, $S^m_{1,0}\subset C^\infty(\Bbb R^n\times\Bbb R^n)$
refers to the Euclidean space symbol class,
defined by the symbol inequalities
\begin{equation}\label{EQ:su2-rn-symbols}
   \left| \partial_\xi^\alpha \partial_x^\beta p(x,\xi) \right|
  \leq C\ \langle\xi\rangle^{m-|\alpha|},
\end{equation} 
for all multi-indices $\alpha,\beta\in\Bbb N_0^n$,
$\N_0=\{0\}\cup\N$,
where $\jp{\xi}=(1+|\xi|^2)^{1/2}$, and where
constant $C$ is independent of $x,\xi\in\Bbb R^n$
but may depend on $\alpha,\beta,p,m$.
On a compact Lie group $G$ we
define the class $\Psi^m(G)$ to be the usual H\"ormander's
class of pseudo-differential operators
of order $m$. Thus, operator
$A$ belongs to $\Psi^m(G)$ if its integral kernel 
$K(x,y)$ is smooth
outside the diagonal $x=y$ and if in (all) local coordinates
operator $A$ is a pseudo-differential operator on $\Rn$
with symbol $p(x,\xi)$ satisfying estimates 
\eqref{EQ:su2-rn-symbols}. We refer to
\cite{Hormander,Hormander2} for the historic development
of this subject.

It is a natural idea to build pseudo-differential operators
out of smooth families of convolution operators on Lie groups.
In this paper, we strive to develop the convolution 
approach
into a symbolic quantization, which always provides a much more
convenient framework for the analysis of operators.
For this, 
our analysis of operators and their symbols is based on the
representation theory of Lie groups.
This leads to the description of the full symbols of 
pseudo-differential 
operators on Lie groups as sequences of matrices of
growing sizes equal to dimensions of representations.
Moreover, the analysis is not confined to neighborhoods of
the neutral element since it does not rely on the 
exponential mapping and its properties.
We also characterize, in terms of the introduced
quantizations, 
standard H\"ormander's classes
$\Psi^m(G)$ on Lie groups. One of the
advantages of the presented approach
is that we obtain a notion
of full (global) symbols compared with only principal symbols
available in the standard theory via localizations.

To illustrate some ideas,
let us now briefly formulate one of the outcomes of this approach
in the case of the 3-dimensional sphere $\S3$.
Before that we note that if we have a closed
simply-connected 3-dimensional manifold $M$,
then by the recently resolved
Poincar\'e conjecture there is a global
diffeomorphism $M\simeq\S3\simeq\SU2$ that turns $M$
into a Lie group with a group structure induced by $\S3$
(or by $\SU2$).  
Thus, we can use the approach developed
in this paper to immediately obtain the corresponding
global quantization of operators on $M$ with respect to
this induced group product. In fact, all the formulae remain
completely the same since the unitary dual of $\SU2$
(or $\S3$ in the quaternionic $\R^4$) is mapped by this
diffeomorphism as well; for an example of this construction 
in the case of $\S3\simeq\SU2$ see
Section \ref{SEC:Pseudos-on-S3}.
The choice of the group structure on $M$ may be not unique
and is not canonical, but after using the machinery that
we develop for $\SU2$, the corresponding quantization 
can be described entirely in terms of $M$, for an
example see
Theorem \ref{THM:s3-main} for $\S3$ and 
Theorem~\ref{THM:simplecharacterization} for $\SU2$.
In this sense, as different quantizations of operators
exist already on $\Rn$ depending on the choice of the
underlying structure (e.g. Kohn--Nirenberg quantization,
Weyl quantizations, etc.), the possibility to choose 
different group products on $M$ resembles this.
In a subsequent paper we will carry out the detailed
analysis of operators 
on homogeneous spaces and on higher dimensional spheres 
$\Sn\simeq {\rm SO}(n+1)/{\rm SO}(n)$ 
viewed as homogeneous spaces. Although we do not have 
general analogues of the diffeomorphic Poincar\'e
conjecture in higher dimensions, this will cover 
cases when $M$ is a convex surface or a surface with
positive curvature tensor, as well as more general
manifolds in terms of their Pontryagin class, etc.

To fix the notation for the Fourier analysis on $\S3$,
let $t^l:\S3\to U(2l+1)\subset \C^{(2l+1)\times(2l+1)}$, 
$l\in\frac12\N_0$, be a family of group homomorphisms,
which are the irreducible continuous 
(and hence smooth) unitary representations
of $\S3$ when it is endowed with the $\SU2$ structure
via the quaternionic product,
see Section~\ref{SEC:Pseudos-on-S3} for details.
The Fourier coefficient $\widehat{f}(l)$
of $f\in C^\infty(\S3)$ is defined by
$
  \widehat{f}(l) = \int_\S3 f(x)\ t^l(x)^*\ {\rm d}x,
$
where the integration is performed with respect to the
Haar measure,
so that  
$\widehat{f}(l)\in\C^{(2l+1)\times (2l+1)}.$
The corresponding Fourier series is given by
\begin{eqnarray*}
  f(x) & = & \sum_{l\in\frac12\N_0} (2l+1)
  \ {\rm Tr}\left( \widehat{f}(l)\ t^l(x) \right).
\end{eqnarray*}
Now, if $A:C^\infty(\S3)\to C^\infty(\S3)$ is
a continuous linear operator,
we define its full symbol as a mapping
$$
  (x,l)\mapsto \sigma_A(x,l),\quad
  \sigma_A(x,l) = t^l(x)^* (A t^l)(x)\in \C^{(2l+1)\times
  (2l+1)}.
$$
Then we have the representation of operator $A$ in the form
$$
  Af(x) = \sum_{l\in\frac12\N_0} (2l+1)
  \ {\rm Tr}\left( t^l(x)\ \sigma_A(x,l)\ \widehat{f}(l)
  \right),
$$
see Theorem \ref{THM:su2-symbol}.
We also note that if
$$
 Af(x)=\int_\S3 K_A(x,y)\ f(y) \drm y
 =\int_\S3 f(y)\ R_A(x,y^{-1}x) \drm y,
$$
where $R_A$ is the right convolution kernel of $A$, then
$\sigma_A(x,l)=\int_\S3 R_A(x,y)\ t^l(y)^* \drm y$
by Theorem \ref{THM:su2-symbol2},
where, as usual, 
the integration is performed with respect to the
Haar measure with a standard distributional interpretation.

One of the arising
fundamental questions is what condition on the matrix
symbols $\sigma_A$ characterize operators from 
H\"ormander's class $\Psi^m(\S3)$. For this, we
introduce symbol class $S^m(\S3)$.
We write $\sigma_A\in S^m(\S3)$ if the corresponding
kernel $K_A(x,y)$ is smooth 
outside the diagonal $x=y$ and if 
we have the estimate
\begin{equation}\label{EQ:S3-symbols}
  \left| \triangle_l^\alpha \partial_x^\beta 
  \sigma_{A_u}(x,l)_{ij} \right|
        \leq C_{A\alpha\beta mN}\  
       (1+|i-j|)^{-N} (1+l)^{m-|\alpha|},
\end{equation} 
for every $N\geq0$, every $u\in\S3$, and all multi-indices
$\alpha,\beta$, 
where symbol $\sigma_{A_u}$ is 
the symbol of operator
$A_u f=A(f\circ\va_u)\circ\va_u^{-1}$, where
$\va_u(x)=xu$ is the quaternionic product.
Symbols of $A_u$ and $A$ can be shown to be related by
formula
$
 \sigma_{A_u}(x,l)=t^l(u)^*\sigma_A(xu^{-1},l)\ t^l(u).
$
We notice that imposing the same conditions
on all symbols $\sigma_{A_u}$ in \eqref{EQ:S3-symbols}
simply refers to the 
well-known fact that
the class $\Psi^m(\S3)$ should be in
particular ``translation"-invariant
(i.e. invariant under the changes of
variables induced by quanternionic products $\va_u$), namely that
$A\in\Psi^m(\S3)$ if and only if $A_u\in\Psi^m(\S3)$,
for all $u\in\S3$.
Condition \eqref{EQ:S3-symbols} is the growth condition
with respect to the
{\em quantum number} 
$l$ combined with a rather striking condition that
matrices $\sigma_A(x,l)$ must have a rapid off-diagonal
decay.
We also write 
$\triangle_l^\alpha=\triangle_+^{\alpha_1} \triangle_-^{\alpha_2}
\triangle_0^{\alpha_3}$, where operators
$\triangle_+, \triangle_-, \triangle_0$ are discrete 
difference operators acting on matrices $\sigma_A(x,l)$
in variable $l$, and explicit formulae for
them and their properties are given in Section~\ref{SEC:su2}.
With this definition, we have the following characterization:
\begin{thm}
\label{THM:s3-main}
We have $A\in\Psi^m(\S3)$
if and only if $\sigma_A\in S^m(\S3)$.
\end{thm} 
The proof of this theorem is based on the detailed analysis of
pseudo-differential operators and their symbols on Lie
group $\SU2$ where we can use its representation theory
and geometric information
to derive the corresponding characterization
of pseudo-differential operators.
We note that this approach works globally on the whole
sphere, since the version of the Fourier analysis is
different from the one in e.g. \cite{She75,St88, She90}
which covers only a hemisphere, with singularities at the
equator.

In our analysis on a Lie group $G$, 
at some point we have to make a choice 
whether to work with
left- or right-convolution kernels. Since
left-invariant operators on $C^\infty(G)$
correspond to right-convolutions $f\mapsto f\ast k$,
once we decide to
identify the Lie algebra $\mathfrak g$ of $G$ with
the left-invariant vector fields on $G$,
it becomes most natural to work with right-convolution kernels
in the sequel, and to define symbols as we do in
Definition~\ref{DEF:su2-symbols-on-G}.

Finally, we mention that the more extensive analysis 
can be carried out in the case of commutative Lie groups.
The main simplification in this case is that full symbols
are just complex-valued scalars (as opposed to 
being matrix-valued in
the non-commutative case) because the continuous irreducible
unitary representations are all one-dimensional. In particular,
we can mention the well-known fact that pseudo-differential 
operators $A\in\Psi^m(\Tn)$ on the $n$-torus
can be globally characterised by conditions 
\begin{equation}\label{EQ:su2-torus-symbols}
   \left| \triangle_\eta^\alpha \partial_x^\beta p(x,\eta) 
   \right|
  \leq C\ (1+|\eta|)^{m-|\alpha|},
\end{equation} 
for all $\eta\in\Z^n$, and all
multi-indices $\alpha,\beta\in\Bbb N_0^n$,
where difference operators
$\triangle^\alpha_\eta=\triangle_{\eta_1}^{\alpha_1}\cdots
\triangle_{\eta_n}^{\alpha_n}$ are defined
by $\triangle_{\eta_j}p(x,\eta)=p(x,\eta+e_j)-p(x,\eta)$,
$(e_j)_k=\delta_{jk}$, for all $1\leq j,k\leq n$, etc. 
If we denote by $S^m(\Tn\times\Zn)$ the class of functions
$p:\Tn\times\Zn\to\C$ satisfying \eqref{EQ:su2-torus-symbols},
then we have $Op S^m(\Tn\times\Zn)=\Psi^m(\Tn)$,
see e.g. 
\cite{Agranovich2,McLean,TurunenVainikko,RT07,RT08} 
with different proofs, as well as
numerical application of this description in e.g.
\cite{SV02,SW87}. We note that in \cite{RT08},
more general symbol classes as well as 
analogues of Fourier integral operators on the torus 
and toroidal microlocal analysis were developed
using the so-called toroidal quantization, which is the 
torus version of the quantization developed here.

It is also known that globally defined symbols of
pseudo-differential operators can be introduced on
manifolds in the presence of a connection which allows
one to use a suitable globally defined phase function,
see e.g. \cite{Wi80,Sa97,Sh05}. However, on a compact
Lie groups the use of the groups structure allows one
to develop a theory parallel to those of $\Rn$ and $\Tn$
in the sense that the Fourier analysis is well adopted
to the underlying representation theory.
Some elements of such theory were discussed in
\cite{Turunen01} and in the PhD thesis of the second
author, and a consistent development 
from different points of view will eventually appear in
\cite{RT-book}.

The global quantization introduced in this paper provides
a relatively easy to use approach to deal with problems
on $\Sn$ (and on more general Lie groups) which depend on
lower order terms of the symbol. Thus, applications to
global hypoellipticity, global solvability and other 
problems in the global setting will appear in the sequel
of this paper.

In this paper, the commutator of matrices 
$X,Y\in\Bbb C^{n\times n}$
will be denoted by
$
  [X,Y]=XY-YX.
$
On $\SU2$, the conventional abbreviations in 
summation indices are

$$
  \sum_{l}  =  \sum_{l\in\frac{1}{2}\Bbb N_0}, \qquad\qquad
  \sum_{l}\sum_{m,n}  =  \sum_{l\in\frac{1}{2}\Bbb N_0}
  \quad\sum_{|m|\leq l,\ l+m\in\Bbb Z}
  \quad\sum_{|n|\leq l,\ l+n\in\Bbb Z},
$$
where $\N_0=\{0\}\cup\N=\{0,1,2,\cdots\}.$
The space of all linear mappings from a finite dimensional
vector space $\Hcal$ to itself will be denoted by
${\rm End}(\Hcal)$. As usual, a mapping $U\in\Lcal(\Hcal)$
is called unitary if $U^*=U^{-1}$ and the 
space of all unitary linear mappings
on a finite dimensional inner product space $\Hcal$
will be denoted by $\Ucal(\Hcal)$.

\section{Full symbols on general compact Lie groups}

Let $G$ be a compact Lie group,
not necessarily just $\SU2$.
Let us endow ${\mathcal D}(G)=C^\infty(G)$
with the usual test function topology.
For a continuous linear operator $A:C^\infty(G)\to C^\infty(G)$,
let $K_A,L_A,R_A\in{\mathcal D}'(G\times G)$
denote respectively the Schwartz, left-convolution and right-convolution
kernels, i.e.
\begin{multline}\label{EQ:SU2-kernels}
  Af(x)
   =  \int_G K_A(x,y)\ f(y)\ {\rm d}y = \\
   =  \int_G L_A(x,xy^{-1})\ f(y)\ {\rm d}y 
   =  \int_G f(y)\ R_A(x,y^{-1}x)\ {\rm d}y
\end{multline}
in the sense of distributions. 
To simplify the notation in the sequel, we will often write
integrals in the sense of distributions, with a standard
distributional interpretation.
Notice that
\begin{equation*}
  R_A(x,y) = L_A(x,xyx^{-1}),
\end{equation*}
and that left-invariant operators on $C^\infty(G)$
correspond to right-convolutions $f\mapsto f\ast k$.
Since we identify the Lie algebra $\mathfrak g$ of $G$ with
the left-invariant vector fields on $G$,
it will be most natural to study right-convolution kernels
in the sequel.

Let us begin with fixing the notation concerning
Fourier series on a compact group $G$
(for general background on the representation theory
we refer to e.g. \cite{HR}).
In the sequel, let $\Rep(G)$ denote the set of all 
strongly continuous irreducible
unitary representations of $G$. 
In this paper, whenever we mention unitary representations
(of a compact Lie group G), we always mean strongly
continuous irreducible unitary representations,
which are then also automatically smooth.
Let $\widehat{G}$ denote the unitary dual of $G$, i.e.
the set of equivalence classes of irreducible unitary 
representations from $\Rep(G)$.
Let $[\xi]\in\widehat{G}$ denote the equivalence class
of an irreducible unitary 
representation $\xi:G\to{\mathcal U}({\mathcal H}_\xi)$;
the representation space ${\mathcal H}_\xi$ is finite-dimensional
since $G$ is compact, and we set $\dim(\xi)=\dim\Hcal_\xi$.
We will always equip compact Lie groups with the 
Haar measure, i.e. the uniquely
determined bi-invariant Borel regular probability measure.
Let us define the Fourier coefficient
$\widehat{f}(\xi) \in {\rm End}({\mathcal H}_\xi)$
of $f\in L^1(G)$ by
\begin{equation}\label{EQ:su2-group-FT}
  \widehat{f}(\xi) := \int_G f(x)\ \xi(x)^\ast\ {\rm d}x;
\end{equation}
more precisely,
$$
  ( \widehat{f}(\xi)u,v )_{{\mathcal H}_\xi}
  = \int_G f(x) \left( \xi(x)^\ast u,v \right)_{{\mathcal H}_\xi} 
  {\rm d}x
  = \int_G f(x) \left( u, \xi(x)v \right)_{{\mathcal H}_\xi} 
  {\rm d}x
$$
for all $u,v\in{\mathcal H}_\xi$, where $(\cdot,\cdot)_{\Hcal_\xi}$
is the inner product of $\Hcal_\xi$. Notice that
$\xi(x)^\ast=\xi(x)^{-1}=\xi(x^{-1})$.

\begin{rem}\label{REM:su2-intertwining}
Let $U\in{\rm Hom}(\eta,\xi)$
be an intertwining isomorphism,
i.e. let $U:{\mathcal H}_\eta\to{\mathcal H}_\xi$
be a bijective unitary linear mapping such that
$U \eta(x)=\xi(x) U$ for every $x\in G$.
Then we have
\begin{equation}\label{EQ:su2-intertwine}
  \widehat{f}(\eta) = U^{-1} \widehat{f}(\xi)\ U
  \in {\rm End}({\mathcal H}_\eta).
\end{equation}
Let us also consider the inner automorphisms
$$
  \phi_u=(x\mapsto u^{-1} x u):G\to G,
$$
where $u\in G$.
If $\xi\in \Rep(G)$ then we also have 
\begin{multline}
\label{EQ:su2-intertwine2}
  \widehat{f\circ\phi_u}(\xi)
   =  \int_G f(u^{-1}xu)\ \xi(x)^\ast\ {\rm d}x 
   =  \int_G f(x)\ \xi(uxu^{-1})^\ast \ {\rm d}x \\
   =  \xi(u) \int_G f(x)\ \xi(x)^\ast\ {\rm d}x\ 
      \xi(u)^\ast 
   =  \xi(u)\ \widehat{f}(\xi)\ \xi(u)^\ast.
\end{multline}
\end{rem}

\begin{rem}\label{REM:su2-convs}
If $f,g\in L^1(G)$ then
\begin{multline*}
  \widehat{f\ast g}(\xi)
   =  \int_G f\ast g(x)\ \xi(x)^\ast\ {\rm d}x 
   =  \int_G \int_G f(xy^{-1}) g(y)\ {\rm d}y\ \xi(x)^\ast
  \ {\rm d}x = \\
   =  \int_G g(y)\ \xi(y)^\ast
  \int_G f(xy^{-1})\ \xi(xy^{-1})^\ast\ {\rm d}x\ {\rm d}y 
   =  \widehat{g}(\xi)\ \widehat{f}(\xi),
\end{multline*}
which in general differs from $\widehat{f}(\xi)\ \widehat{g}(\xi)$.
This order exchange is due to the definition of the 
Fourier coefficients,
where we chose the integration of the function with respect to
$\xi(x)^\ast$ instead of $\xi(x)$.
This choice actually serves us well,
as we chose to identify the Lie algebra $\mathfrak{g}$
with left-invariant vector fields on the Lie group $G$:
namely, a {left}-invariant continuous linear operator
$A:C^\infty(G)\to C^\infty(G)$
can be presented as a {right}-convolution operator
$C_a=(f\mapsto f\ast a)$, resulting in 
convenient expressions like
$$
  \widehat{C_a C_b f\ } =\ \widehat{a}\ \widehat{b}\ \widehat{f}.
$$
\end{rem}

If $\xi:G\to{\rm U}(d)$ is an irreducible 
unitary matrix representation
then
$
  \widehat{f}(\xi) \in\Bbb C^{d\times d} 
$
in \eqref{EQ:su2-group-FT}
has matrix elements
\begin{equation*}
  \widehat{f}(\xi)_{mn} = \int_G f(x)\ 
  \overline{\xi(x)_{nm}}\ {\rm d}x
  \in\Bbb C, \; 1\leq m,n\leq d,
\end{equation*}
where the matrix elements are calculated with respect to the
standard basis of $\C^d$.
If here $f\in L^2(G)$ then
$\widehat{f}(\xi)_{mn} = (f,\xi(x)_{nm})_{L^2(G)}$,
and by the Peter--Weyl Theorem
\begin{equation}\label{EQ:su2-PeterWeyl1}
  f(x)
   = \sum_{[\xi]\in\widehat{G}} {\rm dim}(\xi)
  \ {\rm Tr}\left( \xi(x)\ \widehat{f}(\xi) \right) 
   = \sum_{[\xi]\in\widehat{G}} {\rm dim}(\xi)
  \sum_{m,n=1}^{\dim(\xi)} \xi(x)_{nm}\ \widehat{f}(\xi)_{mn}
\end{equation}
for almost every $x\in G$,
where the summation is understood so that
from each class $[\xi]\in\widehat{G}$
we pick just (any) one representative $\xi\in [\xi]$.
The choice of a representation from the same representation
class is irrelevant due to formula \eqref{EQ:su2-intertwine}
and the presence of the trace in \eqref{EQ:su2-PeterWeyl1}.

\begin{defn}[Symbols of pseudo-differential operators
on $G$]\label{DEF:su2-symbols-on-G}
Let $\xi:G\to{\mathcal U}({\mathcal H}_\xi)$
be an irreducible unitary representation.
The {\it symbol} of a linear continuous
operator $A:C^\infty(G)\to C^\infty(G)$ 
at $x\in G$ and $\xi\in \Rep(G)$
is defined by
$\sigma_A(x,\xi)=\widehat{r_x}(\xi)
\in{\rm End}({\mathcal H}_\xi)$, 
where
$$
  r_x(y) = R_A(x,y)
$$
is the right convolution kernel of $A$ as in 
\eqref{EQ:SU2-kernels}.
Hence
\begin{equation}\label{EQ:su2-right-conv}
  \sigma_A(x,\xi) = \int_G R_A(x,y)\ \xi(y)^\ast\ {\rm d}y
\end{equation}
in the sense of distributions, and operator $A$ can be
represented by its symbol:
\end{defn}

\begin{thm}\label{THM:su2-symbol}
Let the symbol $\sigma_A$ of
a continuous linear operator $A:C^\infty(G)\to C^\infty(G)$
be defined as in Definition~\ref{DEF:su2-symbols-on-G}.
Then
\begin{equation}\label{EQ:su2-psop-def}
  Af(x)
  = \sum_{[\xi]\in\widehat{G}} {\rm dim}(\xi)
  \ {\rm Tr}\left( \xi(x)\ \sigma_A(x,\xi)\ \widehat{f}(\xi)\right).
\end{equation}
for every $f\in C^\infty(G)$ and $x\in G$.
\end{thm}

\begin{proof}
Let us define a right-convolution operator
$A_{x_0}\in{\Lcal}(C^\infty(G))$ by kernel
$R_A(x_0,y) = r_{x_0}(y)$, i.e. by
$$
  A_{x_0}f(x) := \int_G f(y)\ r_{x_0}(y^{-1}x)\ {\rm d}y=
  (f*r_{x_0})(x).
$$
Thus $\sigma_{A_{x_0}}(x,\xi) = \widehat{r_{x_0}}(\xi) 
= \sigma_A(x_0,\xi)$,
so that by \eqref{EQ:su2-PeterWeyl1} we have 
\begin{eqnarray*}
  A_{x_0}f(x)
  & = & \sum_{[\xi]\in\widehat{G}} {\rm dim}(\xi)
  \ {\rm Tr}\left( \xi(x)\ \widehat{A_{x_0}f}(\xi) \right) \\
  & = & \sum_{[\xi]\in\widehat{G}} {\rm dim}(\xi)
  \ {\rm Tr}\left( \xi(x)\ \sigma_A(x_0,\xi)\ \widehat{f}(\xi) \right),
\end{eqnarray*}
where we used that $\widehat{f*r_{x_0}}=\widehat{r_{x_0}}
\widehat{f}$ by Remark~\ref{REM:su2-convs}.
This implies the result, because
$Af(x) = A_x f(x)$.
\end{proof}
For a symbol $\sigma_A$, the corresponding operator 
$A$ defined by 
\eqref{EQ:su2-psop-def} will be also denoted by 
$Op(\sigma_A)$.

Thus, if $\xi:G\to{\rm U}({\rm dim}(\xi))$
are irreducible unitary matrix representations then
\begin{eqnarray*}
  Af(x) & = & \sum_{[\xi]\in\widehat{G}} {\rm dim}(\xi) 
  \sum_{m,n=1}^{\dim(\xi)}
  \xi(x)_{nm} \left( \sum_{k=1}^{\dim(\xi)}
    \sigma_A(x,\xi)_{mk}\ \widehat{f}(\xi)_{kn} \right).
\end{eqnarray*}
Alternatively,
setting $A\xi(x)_{mn} := (A(\xi_{mn}))(x)$, we have
\begin{equation}\label{EQ:su2-symbol-G-2}
  \sigma_A(x,\xi)_{mn}
  = \sum_{k=1}^{\dim(\xi)} \overline{\xi_{km}(x)}\ (A\xi_{kn})(x),
\end{equation}
$1\leq m,n\leq \dim(\xi)$, which follows from the following
theorem:

\begin{thm}
\label{THM:su2-symbol2}
Let the symbol $\sigma_A$ of
a continuous linear operator $A:C^\infty(G)\to C^\infty(G)$
be defined as in Definition~\ref{DEF:su2-symbols-on-G}.
Then
\begin{equation}\label{EQ:su2-symbol-G-1}
  \sigma_A(x,\xi)
  = \xi(x)^\ast (A\xi)(x).
\end{equation}
\end{thm} 
\begin{proof}
Working with representations $\xi:G\to{\rm U}({\rm dim}(\xi))$,
we have
\begin{eqnarray*}
  \sum_{k=1}^{\dim(\xi)} \overline{\xi_{km}(x)}\ (A\xi_{kn})(x)
  &=& \sum_k \overline{\xi_{km}(x)}
  \sum_{[\eta]\in\widehat{G}} \dim(\eta)\ {\rm Tr}\left( \eta(x)
    \ \sigma_A(x,\eta)\ \widehat{\xi_{kn}}(\eta) \right) \\
  &=& \sum_k \overline{\xi_{km}(x)}
  \sum_{[\eta]\in\widehat{G}} \dim(\eta)
  \sum_{i,j,l} \eta(x)_{ij}\ \sigma_A(x,\eta)_{jl}
  \ \widehat{\xi_{kn}}(\eta)_{li} \\
  &=& \sum_{k,j} \overline{\xi_{km}(x)}\ \xi(x)_{kj}\ 
  \sigma_A(x,\xi)_{jn} \\
  &=& \sigma_A(x,\xi)_{mn}, 
\end{eqnarray*}
where if $\eta\in [\xi]$ in the sum, we take $\eta=\xi$,
so that 
$\widehat{\xi_{kn}}(\eta)_{li}=\jp{\xi_{kn},\eta_{il}}_{L^2}$,
which equals $\frac{1}{\dim\xi}$ if $\xi=\eta, k=i$ and
$n=l$, and zero otherwise.
\end{proof}

\begin{rem}\label{REM:su2-symbols}
The symbol of $A\in{\Lcal}(C^\infty(G))$ is a mapping
\begin{equation*}
  \sigma_A:G\times \Rep(G)\to\bigcup_{\xi\in \Rep(G)}{\rm End}
  ({\Hcal}_\xi),
\end{equation*}
where $\sigma_A(x,\xi)\in{\rm End}({\Hcal}_\xi)$
for every $x\in G$ and $\xi\in \Rep(G)$.
However, it can be viewed as a {\em mapping on
the space $G\times\widehat{G}$}. Indeed,
let $\xi,\eta\in \Rep(G)$ be equivalent
via an intertwining isomorphism $U\in{\rm Hom}(\xi,\eta)$:
i.e. such that there exists
a linear unitary bijection
$U:{\Hcal}_\xi\to{\Hcal}_\eta$  
such that $\eta(x)\ U = U\ \xi(x)$ for every $x\in G$,
that is $\eta(x) = U\ \xi(x)\ U^\ast$.
Then by Remark~\ref{REM:su2-intertwining} we have
$\widehat{f}(\eta)  =  U\ \widehat{f}(\xi)\ U^\ast$, and
hence also
$$
  \sigma_A(x,\eta)  =  U\ \sigma_A(x,\xi)\ U^\ast.
$$
Therefore, taking any representation from the same 
class $[\xi]\in\widehat{G}$ leads to the same
operator $A$ in view of the trace in 
formula \eqref{EQ:su2-psop-def}.
In this sense we may think that symbol $\sigma_A$
is defined on $G\times\widehat{G}$ instead of $G\times \Rep(G)$.
\end{rem}

Notice that if $A=(f\mapsto f\ast a)$
then $R_A(x,y)=a(y)$ and
$$
  \sigma_A(x,\xi) = \widehat{a}(\xi),
$$
i.e. $\widehat{Af}(\xi) = \widehat{a}(\xi)\ \widehat{f}(\xi)$.
Moreover, if $B=(f\mapsto b\ast f)$
then $L_B(x,y)=b(y)$, $R_B(x,y)=L_B(x,xyx^{-1})=b(xyx^{-1})$,
and by \eqref{EQ:su2-intertwine2} we have
$$
  \sigma_B(x,\xi) = \xi(x)^\ast\ \widehat{b}(\xi)\ \xi(x).
$$

\begin{rem}
Let $\mathfrak g$ be the Lie algebra of a compact Lie group $G$,
and let $n={\rm dim}(G)={\rm dim}(\mathfrak{g})$.
By the exponential mapping $\exp:\mathfrak{g}\to G$,
a neighbourhood of the neutral element $e\in G$
can be identified with a neighbourhood of $0\in\mathfrak g$.
Let ${\mathcal X}^m = S^m_{1\#}\subset S^m_{1,0}$ consist
of the $x$-invariant symbols $(x,\xi)\mapsto p(\xi)$ in $S^m_{1,0}$
with the usual Fr\'echet space topology.
A distribution $k\in{\mathcal D}'(G)$
with a sufficiently small support
is said to belong to space $\widehat{{\mathcal X}^m}$ if
$
  {\rm sing\ supp}(k)\subset\{e\}$ and 
  $\widehat{k}\in{\mathcal X}^m\subset C^\infty(\mathfrak{g}'),
$
where the Fourier transform $\widehat{k}$ is the 
usual Fourier transform
on $\mathfrak{g}\cong\Bbb R^n$, and the dual space
satisfies $\mathfrak{g}'\cong\Bbb R^n$
(and we are using the exponential coordinates for $k(y)$
when $y\approx e\in G$).
If $k\in\widehat{{\mathcal X}^m}$ then the convolution operator
$$
  u\mapsto k\ast u,\quad
  k\ast u(x) = \int_G k(xy^{-1})\ u(y)\ {\rm d}y,
$$
is said to belong to space $OP{\mathcal X}^m$,
which is endowed with the natural Fr\'echet 
space structure obtained from
${\mathcal X}^m$.
Formally, let $k(x,y)=k_x(y)$ be the left-convolution kernel
of a linear operator ${\mathcal K}: C^\infty(G)\to C^\infty(G)$,
i.e.
$$
  {\mathcal K}u(x) = \int_G k_x(xy^{-1})\ u(y)\ {\rm d}y.
$$
In \cite{Taylor84}, M. E. Taylor showed that
$\mathcal K\in \Psi^m(G)$ if and only if the mapping
$$
  \left(x\mapsto (u\mapsto k_x\ast u)\right):G\to OP{\mathcal X}^m
$$
is smooth; here naturally $u\mapsto k_x\ast u$
must belong to $OP{\mathcal X}^m$ for each $x\in G$.
\end{rem} 

In the sequel, we will need conjugation properties of
symbols which we will now analyse for this purpose.
\begin{defn}
Let $\phi:G\to G$ be a diffeomorphism,
$f\in C^\infty(G)$, $A:C^\infty(G)\to C^\infty(G)$ 
continuous and linear.
Then the {\it $\phi$-pushforwards} $f_\phi\in C^\infty(G)$
and $A_\phi:C^\infty(G)\to C^\infty(G)$
are defined by
\begin{eqnarray*}
  f_\phi &:=& f\circ\phi^{-1}, \\
  A_\phi f & := & \left(A(f_{\phi^{-1}})\right)_{\phi}
  \quad = \quad A(f\circ\phi)\circ\phi^{-1}.
\end{eqnarray*}
Notice that
$$
  A_{\phi\circ\psi} = \left(A_\psi\right)_\phi.
$$
From the local theory of pseudo-differential operators,
it is well-known that $A\in\Psi^\mu(G)$
if and only if $A_\phi\in\Psi^\mu(G)$.
\end{defn}

\begin{defn}
For $u\in G$, let $u_L,u_R:G\to G$ be defined by
$$
  u_L(x) := ux\quad {\rm and}\quad u_R(x):= xu.
$$
Then $(u_L)^{-1} = (u^{-1})_L$ and $(u_R)^{-1} = (u^{-1})_R$.
The inner automorphism $\phi_u:G\to G$ defined 
in Remark~\ref{REM:su2-intertwining}
by
$
  \phi_u(x) := u^{-1}xu
$
satisfies $\phi_u=u_L^{-1}\circ u_R = u_R\circ u_L^{-1}$.
\end{defn}

\begin{prop}\label{PROP:conjugatedrightsymbols}
Let $u\in G$, $B=A_{u_L}$, $C=A_{u_R}$ and $F=A_{\phi_u}$.
Then we have the following relations between symbols:
\begin{eqnarray*}
  \sigma_B(x,\xi) & = & \sigma_A(u^{-1}x,\xi), \\
  \sigma_C(x,\xi) & = & \xi(u)^\ast\ 
  \sigma_A(xu^{-1},\xi)\ \xi(u), \\
  \sigma_F(x,\xi) & = & \xi(u)^\ast\ \sigma_A(u x u^{-1},\xi)\ \xi(u).
\end{eqnarray*}
Especially, if $A=(f\mapsto f\ast a)$,
i.e. $\sigma_A(x,\xi)=\widehat{a}(\xi)$,
then
\begin{eqnarray*}
  \sigma_B(x,\xi) & = & \widehat{a}(\xi), \\
  \sigma_C(x,\xi) & = & \xi(u)^\ast\ \widehat{a}(\xi)\ \xi(u) 
   =  \sigma_F(x,\xi).
\end{eqnarray*}
\end{prop}

\begin{proof}
We notice that $F=C_{(u^{-1})_L}$, so it suffices
to consider only operators $B$ and $C$.
For operator $B=A_{u_L}$, we get
\begin{multline*}
  \int_G f(z)\ R_B(x,z^{-1}x)\ {\rm d}z
   =  Bf(x) 
   =  A(f\circ u_L)(u_L^{-1}(x)) = \\
   =  \int_G f(uy)\ R_A(u^{-1}x,y^{-1}u^{-1}x)\ {\rm d}y 
   =  \int_G f(z)\ R_A(u^{-1}x,z^{-1} x)\ {\rm d}z,
\end{multline*}
so $R_B(x,y) = R_A(u^{-1}x,y)$, yielding
$
  \sigma_B(x,\xi)  =  \sigma_A(u^{-1}x,\xi).
$
For operator $C=A_{u_R}$, we get
similarly
$R_C(x,y) = R_A(xu^{-1},uyu^{-1})$, yielding the result.
\end{proof}

Let us finally record how push-forwards by translation affect vector
fields.
\begin{lem}\label{LEM:conjugateddiff}
Let $u\in G$, $Y\in\mathfrak{g}$ and let 
$E=D_Y:C^\infty(G)\to C^\infty(G)$
be defined by $D_Y f (x)
  = \left. \frac{\rm d}{{\rm d}t} f(x\ \exp(tY)) \right|_{t=0}.$
Then
$$
  E_{u_R} = E_{\phi_u} = D_{u^{-1}Yu},
$$
i.e.
$
  D_Y(f\circ u_R)(xu^{-1}) = D_Y(f\circ\phi_u)(uxu^{-1})
  = D_{u^{-1}Yu}f(x).
$
\end{lem}

\begin{proof} We have
\begin{multline*}
  E_{u_R} f(x) =  E(f\circ u_R)(xu^{-1}) 
   =  \frac{\rm d}{{\rm d}t} \left.(f\circ u_R)(xu^{-1}\exp(tY))
    \right|_{t=0} = \\
   =  \frac{\rm d}{{\rm d}t} \left. f(xu^{-1}\exp(tY)u)) 
   \right|_{t=0} 
   =  \frac{\rm d}{{\rm d}t} \left. f(x\exp(tu^{-1}Yu) 
   \right|_{t=0} 
   =  D_{u^{-1}Yu}f(x).
\end{multline*}
Due to the left-invariance, $E_{u_L}=E$, so that
$E_{\phi_u}=(E_{u_L^{-1}})_{u_R}=E_{u_R}=D_{u^{-1}Yu}.$
\end{proof}

\section{Boundedness of pseudo-differential operators
on $L^2(G)$ and $H^s(G)$}

In this section we will state some natural conditions
on the symbol of an operator $A:C^\infty(G)\to C^\infty(G)$
to guarantee the boundedness on Sobolev spaces.
The Sobolev space $H^s(G)$ of order $s\in\Bbb R$
can be defined via a smooth partition of unity
of the closed manifold $G$.

The {Hilbert--Schmidt inner product}
of $A,B\in\mathbb C^{m\times n}$ is
\begin{equation*}
  \langle A,B\rangle_{HS} := {\rm Tr}(B^\ast A)
  = \sum_{i=1}^m \sum_{j=1}^n \overline{B_{ij}} A_{ij},
\end{equation*}
with the corresponding norm
$\|A\|_{HS} := \langle A,A\rangle_{HS}^{1/2}$,
and the {operator norm} 
\begin{equation*}
  \|A\|_{op} := \sup \left\{ \|Ax\|_{HS}:
    \ x\in\mathbb C^{n\times 1},\ \|x\|_{HS} \leq 1 \right\}
  = \|A\|_{\ell^2\to\ell^2}.
\end{equation*}
Let $A,B\in\mathbb C^{n\times n}$. Then we have
$
  \|AB\|_{HS} \leq
  \|A\|_{op}\ \|B\|_{HS}.
$
Moreover, we also have
$
  \|A\|_{op}=
  \sup\left\{ \|AX\|_{HS}:
    \ X\in\mathbb C^{n\times n},\ \|X\|_{HS}\leq 1 \right\}.
$
By this, taking the Fourier transform and using 
Plancherel's formula (see e.g. \cite{Segal1}), we get
\begin{equation}\label{EQ:su2-l2-norm-mult}
  \| g\mapsto f\ast g\|_{{\mathcal L}(L^2(G))} = 
  \| g\mapsto g\ast f\|_{{\mathcal L}(L^2(G))} =
  \sup_{\xi\in \Rep(G)} \|\widehat{f}(\xi)\|_{op},
\end{equation}
by Remark~\ref{REM:su2-convs}.
We also note that 
$\|\widehat{f}(\xi)\|_{op}=\|\widehat{f}(\eta)\|_{op}$
if $[\xi]=[\eta]\in\widehat{G}$.

Let us first consider a condition on
the symbol for the corresponding operator
to be bounded on $L^2(G)$.

\begin{thm}\label{THM:su2-L2}
Let $G$ be a compact Lie group of dimension $n$ and
let $k$ be an integer such that $k>n/2$.
Let $A$ be an operator with symbol $\sigma_A$ 
defined as in Definition~\ref{DEF:su2-symbols-on-G}.
Assume that
there is a constant $C$ such that
$$\|\partial_x^\alpha\sigma_A(x,\xi)\|_{op}\leq C$$
for all $x\in G$, all $\xi\in\Rep(G),$ and all 
$|\alpha|\leq k$, where
$\partial_x^\alpha=\partial_1^{\alpha_1}\cdots 
\partial_n^{\alpha_n}$, and $\partial_1,\ldots,\partial_n$
are first-order differential operators
corresponding to a basis of the Lie algebra of $G$.
Then $A$ is bounded from $L^2(G)$ to $L^2(G)$.
\end{thm}

\begin{proof}
Let $Af(x)=(f*r_A(x))(x)$, where $r_A(x)(y)=R_A(x,y)$
is the right-convolution kernel of $A$.
Let $A_y f(x)=(f*r_A(y))(x)$, so that $A_x f(x)=Af(x)$.
Then
$$
  \| Af \|_{L^2(G)}^2  = 
        \int_G |A_x f(x)|^2\ {\rm d}x \\
  \leq  \int_G \sup_{y\in G} |A_y f(x)|^2\ {\rm d}x,
$$
and by an application of the Sobolev embedding theorem we get
$$
  \sup_{y\in G} |A_y f(x)|^2 \leq
  C \sum_{|\alpha|\leq k} \int_G 
  |\partial_y^\alpha A_y f(x)|^2
        \ {\rm d}y.
$$
Therefore, using the
Fubini theorem to change the order of integration, we obtain
\begin{eqnarray*}
    \| Af \|_{L^2(G)}^2
  &\leq & C \sum_{|\alpha|\leq k} \int_G \int_G
        | \partial_y^\alpha A_y f(x) |^2
        \ {\rm d}x\ {\rm d}y \\
  &\leq & C \sum_{|\alpha|\leq k} \sup_{y\in G} \int_G
        | \partial_y^\alpha A_y f(x) |^2\ {\rm d}x \\
  & = & C \sum_{|\alpha|\leq k} \sup_{y\in G}
        \| \partial_y^\alpha A_y f \|_{L^2(G)}^2 \\
  &\leq & C \sum_{|\alpha|\leq k} \sup_{y\in G}
        \| f\mapsto f*\partial_y^\alpha r_A(y)\|_
        {{\mathcal L}(L^2(G))}^2
        \|f\|_{L^2(G)}^2 \\
  &\leq & C \sum_{|\alpha|\leq k} \sup_{y\in G}
   \sup_{[\xi]\in\widehat{G}}
        \| \partial_y^\alpha
         \sigma_A(y,\xi)\|_{op}^2
          \|f\|_{L^2(G)}^2,
\end{eqnarray*}
where the last inequality holds due to 
\eqref{EQ:su2-l2-norm-mult}.
This completes the proof.
\end{proof}

Let $\Lap$ be the bi-invariant 
Laplacian of $G$,
i.e. the Laplace-Beltrami operator corresponding to
the unique (up to scaling) bi-invariant Riemannian metric of $G$.
The Laplacian is symmetric and $I-\Lap$ is positive.
Denote $\Xi = (I-\Lap)^{1/2}$.
Then $\Xi^s\in {\mathcal L}(C^\infty(G))$ and
$\Xi^s\in {\mathcal L}({\mathcal D}'(G))$ 
for every $s\in\mathbb R$.
Let us define
$$
  (f,g)_{H^s(G)} = (\Xi^{s} f, \Xi^{s}g )_{L^2(G)}
  \ \ (f,g\in C^\infty(G)).
$$
The completion of $C^\infty(G)$ with respect to the norm
$f\mapsto \| f\|_{H^s(G)} = (f,f)_{H^s(G)}^{1/2}$
gives us Sobolev space $H^s(G)$ of order $s\in\mathbb R$,
which coincides with the Sobolev space
obtained using any smooth partition of unity
on the compact manifold $G$.
Operator $\Xi^r$ is a Sobolev space isomorphism
$H^s(G)\to H^{s-r}(G)$ for every $r,s\in\mathbb R$.
To formulate the corresponding boundedness result in
Sobolev spaces, let us introduce some notation.

Let $\xi\in \Rep(G)$.
Given $v,w\in \Hcal_\xi$, the function $\xi^{vw}:G\to\Bbb C$
defined by
$$
  \xi^{vw}(x) := \langle \xi(x)v,w\rangle_{\Hcal_\xi}
$$
is not only continuous but even $C^\infty$-smooth.
Let ${\rm span}(\xi)$ denote the linear span of
$
  \left\{ \xi^{vw}:\ v,w\in \Hcal_\xi \right\}.
$
If $\xi\sim\eta$ then ${\rm span}(\xi)={\rm span}(\eta)$;
consequently, we may write
$$
  {\rm span}[\xi] := {\rm span}(\xi) \subset C^\infty(G).
$$
It follows that $-\Lap\xi^{vw}(x)=\lambda_{[\xi]}
\xi^{vw}(x)$, where $\lambda_{[\xi]}\geq 0$, and we denote 
\begin{equation}\label{EQ:su2-jp-bracket}
\jp{\xi}=(1+\lambda_{[\xi]})^{1/2}.
\end{equation}
We note that $\sigma_\Lap(x,\xi)=-\lambda_{[\xi]}I_{\dim\xi},$
where $I_{\dim\xi}$ is the identity mapping on 
$\Hcal_\xi$.

Now we can formulate the main result on Sobolev space boundedness:

\begin{thm}\label{THM:su2-Sobolev}
Let $G$ be a compact Lie group of dimension $n$.
Let $A$ be an operator with symbol $\sigma_A$ 
defined as in Definition~\ref{DEF:su2-symbols-on-G}.
Assume that
there are constants $\mu,C_\alpha\in\R$ such that
$$
  \|\partial_x^\alpha\sigma_A(x,\xi)\|_{op}\leq C_\alpha\ \jp{\xi}^\mu
$$
holds for all $x\in G$, $\xi\in\Rep(G)$ and 
all multi-indices $\alpha$,
where
$\partial_x^\alpha=\partial_1^{\alpha_1}\cdots\partial_n^{\alpha_n}$
is as in Theorem~\ref{THM:su2-L2}.
Then
$A$ is bounded from $H^s(G)$ to $H^{s-\mu}(G)$,
for all $s\in\R$.
\end{thm}

\begin{rem}\label{REM:SobolevL2boundedness}
We shall prove this theorem later in Section~\ref{SEC:psdosymbols},
after introducing tools for symbolic calculus.
However, notice that we may easily obtain a special 
case of this result with $s=\mu$.
Namely, if $\sigma_A$ is as in Theorem~\ref{THM:su2-Sobolev},
then
$$
  \left\| \partial_x^\alpha \left(
      \sigma_A(x,\xi) \langle\xi\rangle^{-\mu}\right) \right\|_{op}
  \leq C_\alpha
$$
for every multi-index $\alpha$. Here
$\sigma_A(x,\xi)\langle\xi\rangle^{-\mu}=
\sigma_{A\circ\Xi^{-\mu}}(x,\xi)$,
and thus Theorem~\ref{THM:su2-L2} implies that $A\circ\Xi^{-\mu}$
is bounded on $L^2(G)$,
so that $A\in{\mathcal L}(H^\mu(G),L^2(G))$.
\end{rem}

%
%

\section{Preliminaries on $\SU2$}
\label{SEC:preliminaries}

We study the compact group $\SU2$ defined by
$$
  \SU2 = \left\{ u\in\mathbb C^{2\times 2}:
    \ {\rm det}(u) = 1\ {\rm and}\ u^\ast u = I \right\},
$$
where $e=I=\begin{pmatrix} 1 & 0 \\ 0 & 1 
\end{pmatrix}\in\mathbb C^{2\times 2}$
is the identity matrix.
Matrix $u\in\mathbb C^{2\times 2}$ belongs to $\SU2$
if and only if it is of the form
$
  u = \begin{pmatrix}
    \alpha & \beta \\
    -\overline{\beta} & \overline{\alpha}
    \end{pmatrix},
$
where $|\alpha|^2 + |\beta|^2 =1$.
We will now fix the notation concerning the representations
of $\SU2$.
Let us identify $z=(z_1,z_2)\in\mathbb C^2$ with matrix
$z=\begin{pmatrix} z_1 & z_2 \end{pmatrix}\in\mathbb C^{1\times 2}$,
and let $\mathbb C[z_1,z_2]$ be the space
of two-variable polynomials 
$f:\mathbb C^2\to\mathbb C$.
Consider mappings
\begin{equation*}
  t^l:\SU2\to {\rm GL}(V_l),\quad (t^l(u) f)(z) = f(zu),
\end{equation*}
where $l\in\frac{1}{2}\mathbb N_0$ may be called the
quantum number, and where
$V_l$ is the ($2l+1$)-dimensional
subspace of $\mathbb C[z_1,z_2]$
containing the homogeneous polynomials
of order $2l\in\mathbb N_0$, i.e.
\begin{equation*}
  V_l = \left\{ f\in\mathbb C[z_1,z_2]:
    \ f(z_1,z_2)
    =\sum_{k=0}^{2l} a_k z_1^k z_2^{2l-k},
    \quad \{a_k\}_{k=0}^{2l}\subset\mathbb C \right\}.
\end{equation*}
Then the family $\{t^l\}_{l\in\frac12\N_0}$ is 
the family of irreducible unitary representations
of $\SU2$ such that any other 
irreducible unitary representation of $\SU2$
is equivalent to one of $t^l$.
The collection
$\left\{q_{lk}:\ k\in\{-l,-l+1,\cdots,+l-1,+l\} \right\}$
is a basis for the representation space $V_l$,
where
$$
  q_{lk}(z) = \frac{z_l^{l-k} z_2^{l+k}}{\sqrt{(l-k)!(l+k)!}}.
$$
Let us give the matrix elements $t^l_{mn}(u)$
of $t^l(u)$ with respect to this basis, where 
\eqref{EQ:tlmn-in-derivatives} is well-known and
\eqref{EQ:tlmn-explicit} follows from it.

\begin{prop}\label{matrixelements}
Let
$
  u = \begin{pmatrix} a & b 
  \\ c & d \end{pmatrix}.
$
Then
\begin{equation}\label{EQ:tlmn-in-derivatives}
 t^l_{mn} (u)
  =
  \left(\frac{{\rm d}}{{\rm d}z_1}\right)^{l-m}
  \left(\frac{{\rm d}}{{\rm d}z_2}\right)^{l+m}
  \frac{
  (z_1 a + z_2 c)^{l-n} (z_1 b + z_2 d)^{l+n}}
  {\sqrt{(l-m)!(l+m)!(l-n)!(l+n)!}},
\end{equation}
where
$$
  P^l_{mn}(x) = c^l_{mn}
  \frac{(1-x)^{(n-m)/2}}{(1+x)^{(m+n)/2}}
  \ \left(\frac{{\rm d}}{{\rm d}x}\right)^{l-m}
  \left[ (1-x)^{l-n} (1+x)^{l+n} \right]
$$
with
$$
  c^l_{mn} = 2^{-l} \frac{(-1)^{l-n}\ 
  {\rm i}^{n-m}}{\sqrt{(l-n)!\ (l+n)!}}
  \ \sqrt{\frac{(l+m)!}{(l-m)!}}.
$$
Moreover, we have
\begin{multline}\label{EQ:tlmn-explicit}
 t^l_{mn}(u) =
 \sqrt{\frac{(l-m)!(l+m)!}{(l-n)!(l+n)!}}\times \\
 \times\sum_{i=\max\{0,n-m\}}^{\min\{l-n,l-m\}}
 \frac{(l-n)!(l+n)!}{i! (l-n-i)!(l-m-i)!(n+m+i)!}
 a^i b^{l-m-i} c^{l-n-i} d^{n+m+i}.
\end{multline}
\end{prop}

On a compact group $G$,
a function $f:G\to\mathbb C$ is called a trigonometric polynomial
if its translates span the finite-dimensional vector space,
i.e. if
$$
  \dim{\rm span} 
  \left\{ (x\mapsto f(y^{-1}x)):G\to\mathbb C\ |\ y\in G \right\}
  < \infty.
$$
A trigonometric polynomial can be expressed as a linear combination
of matrix elements of irreducible unitary representations.
Thus a trigonometric polynomial is continuous,
and on a Lie group even $C^\infty$-smooth.
Moreover,
trigonometric polynomials form an algebra
with the usual pointwise multiplication.
On $\SU2$, actually,
$$
  t^{l'}_{m'n'}\ t^l_{mn} = \sum_{k=|l-l'|}^{l+l'} 
  C^{l l' (l+k)}_{m' m (m'+m)}
        \ C^{l l' (l+k)}_{n' n (n'+n)}
        \ t^{l+k}_{(m'+m)(n'+n)},
$$
where $C^{l l' (l+k)}_{m' m (m'+m)}$
are Clebsch-Gordan coefficients, for which there are
explicit formulae, see e.g. \cite{Vilenkin}.
Now we are going to give basic multiplication formulae
for trigonometric polynomials $t^l_{mn}:\SU2\to\mathbb C$;
for general multiplication of trigonometric polynomials,
one can use these formulae iteratively.

\begin{thm}\label{multiplication}
Let
$$
\begin{pmatrix} t_{--} & t_{-+} \\ t_{+-} & t_{++} \end{pmatrix}
\equiv t^{1/2} = \begin{pmatrix}
  t^{1/2}_{-1/2,-1/2} & t^{1/2}_{-1/2,+1/2} \\
  t^{1/2}_{+1/2,-1/2} & t^{1/2}_{+1/2,+1/2} \end{pmatrix}
$$
and denote $x^{\pm} := x\pm 1/2$ for $x\in\mathbb R$. Then
\begin{eqnarray*}
  t^l_{mn} t_{--}
  & = &
    \frac{\sqrt{(l-m+1)(l-n+1)}}{2l+1}\ t^{l^+}_{m^- n^-}
  + \frac{\sqrt{(l+m)(l+n)}}{2l+1}\ t^{l^-}_{m^- n^-}, \\
  t^l_{mn} t_{++}
  & = &
    \frac{\sqrt{(l+m+1)(l+n+1)}}{2l+1}\ t^{l^+}_{m^+ n^+}
  + \frac{\sqrt{(l-m)(l-n)}}{2l+1}\ t^{l^-}_{m^+ n^+}, \\
  t^l_{mn} t_{-+}
  & = &
    \frac{\sqrt{(l-m+1)(l+n+1)}}{2l+1}\ t^{l^+}_{m^- n^+}
  - \frac{\sqrt{(l+m)(l-n)}}{2l+1}\ t^{l^-}_{m^- n^+}, \\
  t^l_{mn} t_{+-}
  & = &
    \frac{\sqrt{(l+m+1)(l-n+1)}}{2l+1}\ t^{l^+}_{m^+ n^-}
  - \frac{\sqrt{(l-m)(l+n)}}{2l+1}\ t^{l^-}_{m^+ n^-}.
\end{eqnarray*}
\end{thm}
These formulae imply, in particular, that expressions similar
to these will appear naturally in the developed quantization
of operators on $\SU2$.

\section{Left-invariant differential operators on $\SU2$}

Let us analyse first-order partial differential operators on $\SU2$
from the point of view of pseudo-differential operators and their
global quantization.
Homomorphisms $\omega:\Bbb R\to \SU2$ are called 
{one-parametric subgroups},
and they are of the form $\omega=(t\mapsto\exp(tY))$
for $Y=\omega'(0)\in\mathfrak{su}(2)$.
As usual, we identify the Lie algebra $\mathfrak{su}(2)$
with the left-invariant vector fields on $\SU2$,
by associating $Y\in\mathfrak{su}(2)$
to the left-invariant operator
$D_Y:C^\infty(\SU2)\to C^\infty(\SU2)$
defined by
\begin{equation}\label{diffopLie}
  D_Y f (x)
  = \left. \frac{\rm d}{{\rm d}t} f(x\ \exp(tY)) \right|_{t=0}.
\end{equation}
\begin{rem}\label{REM:su2-symbols-vfs}
Notice first that vector field ${\rm i}D_Y$ is symmetric
on an arbitrary $G$:
\begin{multline*}
  \left( {\rm i}D_Y f,g\right)_{L^2(G)}
   =  \int_G ({\rm i}D_Y f)(x)\ \overline{g(x)}\ {\rm d}x 
   =  -{\rm i} \int_G f(x)\ \overline{D_Y g(x)}\ {\rm d}x 
   =  \left( f,{\rm i}D_Y g \right)_{L^2(G)}.
\end{multline*}
Hence it is always possible to choose a representative 
$\xi\in \Rep(G)$
from each $[\xi]\in\widehat{G}$ such that
$\sigma_{{\rm i}D_Y}(x,\xi)$
is a diagonal matrix
$\begin{pmatrix}\lambda_1 & & \\ 
& \ddots & \\ & & \lambda_{{\rm dim}(\xi)}
\end{pmatrix}$,
with diagonal entries $\lambda_j\in\Bbb R$,
which follows because symmetric matrices can be
diagonalised by unitary matrices.
Notice that then also
$
  [\sigma_{{\rm i}D_Y},\sigma_A](x,\xi)_{mn}
  = (\lambda_m-\lambda_n)\ \sigma_A(x,\xi)_{mn}.
$
\end{rem}
In the case of $\SU2$, we will simplify the notation
writing $\widehat{f}(l)$ instead of $\widehat{f}(t^l)$, etc.,
since we can take a representative $t^l$ in each 
equivalence class in $\widehat{\SU2}$.
\begin{defn}
Let us define one-parametric subgroups
$\omega_1,\omega_2,\omega_3:\mathbb R\to \SU2$ by
\begin{eqnarray*}
  \omega_1(t) & = &
        \begin{pmatrix}
                \cos(t/2) & {\rm i}\sin(t/2) \\
                {\rm i}\sin(t/2) & \cos(t/2)
        \end{pmatrix}, \\
  \omega_2(t) & = &
        \begin{pmatrix}
                \cos(t/2) & -\sin(t/2) \\
                \sin(t/2) & \cos(t/2)
        \end{pmatrix}, \\
  \omega_3(t) & = &
        \begin{pmatrix}
                {\rm e}^{{\rm i}t/2} & 0 \\
                0 & {\rm e}^{-{\rm i}t/2}
        \end{pmatrix}.
\end{eqnarray*}
Let $Y_j := \omega_j'(0)$, i.e.
\begin{equation*}
  Y_1 = \frac{1}{2}
        \begin{pmatrix}
                0 & {\rm i} \\
                {\rm i} & 0
        \end{pmatrix}, \quad
  Y_2 = \frac{1}{2}
        \begin{pmatrix}
                0 & -1 \\
                1 & 0
        \end{pmatrix}, \quad
  Y_3 = \frac{1}{2}
        \begin{pmatrix}
                {\rm i} & 0 \\
                0 & -{\rm i}
        \end{pmatrix}.
\end{equation*}
Matrices $Y_1,Y_2,Y_3$ constitute a basis for
the real vector space $\mathfrak{su}(2)$.
Notice that
$$
  [Y_1,Y_2]=Y_3,\quad [Y_2,Y_3]=Y_1,\quad [Y_3,Y_1]=Y_2.
$$
Let us define differential operators $D_j := D_{Y_j}$.
\end{defn}

We note that matrices $\frac{2}{\irm}Y_j$, $j=1,2,3$, are known
as Pauli (spin) matrices in physics. 
It can be also noted
that $\mathfrak{k}={\rm span}\{Y_3\}$ and
$\mathfrak{p}={\rm span}\{Y_1,Y_2\}$ form a Cartan
pair of the Lie algebra $\mathfrak{su}(2)$.

\begin{prop}\label{PROP:1parametersubgroupconjugation}
Let $w_j=\omega_j(\pi/2)$ and $t\in\Bbb R$.
Then
$$
  w_1\ \omega_2(t)\ w_1^{-1}  =  \omega_3(t),\;
  w_2\ \omega_3(t)\ w_2^{-1}  =  \omega_1(t),\;
  w_3\ \omega_1(t)\ w_3^{-1}  =  \omega_2(t).
$$
The differential versions of these formulae are
\begin{equation*}
  w_1\ Y_2\ w_1^{-1}  =  Y_3,\;
  w_2\ Y_3\ w_2^{-1}  =  Y_1,\;
  w_3\ Y_1\ w_3^{-1}  =  Y_2.
\end{equation*}
\end{prop}
The proof is straightforward and follows simply by 
multiplying these matrices.

\begin{prop}\label{PROP:SU2-vector-fields}
We have
\begin{equation*}
  (D_3)_{(w_1)_R}  =  D_2,\
  (D_1)_{(w_2)_R}  =  D_3,\
  (D_2)_{(w_3)_R}  =  D_1.
\end{equation*}
Symbols of operators $D_1,D_2$ can be turned to that of $D_3$
by taking suitable conjugations:
\begin{eqnarray}
  \sigma_{D_1}(x,l) & = & t^l(w_2)\ \sigma_{D_3}(x,l)\ t^l(w_2)^\ast,
  \label{linkD1D3} \\
  \sigma_{D_2}(x,l) & = & t^l(w_1)^\ast\ \sigma_{D_3}(x,l)\ t^l(w_1).
  \label{linkD2D3}
\end{eqnarray}
Moreover, if $D\in \mathfrak{su}(2)$ there is $u\in\SU2$ such that
$\sigma_D(l) = t^l(u)^\ast\ \sigma_{D_3}(l)\ t^l(u)$.
\end{prop}
\begin{proof}
Combining Lemma~\ref{LEM:conjugateddiff}
with Proposition~\ref{PROP:1parametersubgroupconjugation},
we see that 
$
  (D_3)_{(w_1)_R}  =  D_2, $
$  (D_1)_{(w_2)_R}  =  D_3 $ and
$  (D_2)_{(w_3)_R}  =  D_1.$
Since $D_1,D_2,D_3$ are left-invariant operators,
their symbols $\sigma_{D_j}(x,l)$ do not depend on $x\in G$,
and by Proposition~\ref{PROP:conjugatedrightsymbols}
we obtain \eqref{linkD1D3} and \eqref{linkD2D3}.
The last statement follows from
Proposition~\ref{PROP:conjugatedrightsymbols}
since $D$ is a rotation of $D_3$.
\end{proof}

Although operators $D_j$ have meaning as derivatives
with respect to $\frac{\irm}{2}$
Pauli matrices, it will be technically
simpler for us to work with their linear combinations
(see Remark~\ref{REM:su2-creation}, also for
the explanation of the terminology), which we will now
define.
\begin{defn}
Let us define left-invariant first-order partial differential operators
$\partial_+,\partial_-,\partial_0:C^\infty(\SU2)\to 
C^\infty(\SU2)$, called creation, annihilation, and neutral
operators, respectively, by
$$
  \begin{cases}
    \partial_+ := {\rm i} D_1 - D_2, \\
    \partial_- := {\rm i} D_1 + D_2, \\
    \partial_0 := {\rm i} D_3,
    \end{cases}
    \quad {\rm i.e.}\quad
  \begin{cases}
    D_1 = \frac{-{\rm i}}{2} 
    \left(\partial_- + \partial_+ \right) \\
    D_2 = \frac{1}{2} \left( \partial_- - \partial_+ \right), \\
    D_3 = -{\rm i} \partial_0.
  \end{cases}
$$
\end{defn}
\begin{rem}
The Laplacian $\Lap$ satisfies $\Lap = D_1^2 + D_2^2 + D_3^2$
and
$
  [\Lap,D_j] = 0
$
for every $j\in\{1,2,3\}$.
Notice that it can be expressed as
$\Lap = - \partial_0^2 - 
(\partial_+\partial_-+\partial_-\partial_+)/2$.
Operators $\partial_+,\partial_-,\partial_0$ satisfy
$[\partial_0,\partial_+]=\partial_+,
[\partial_-,\partial_0]=\partial_-,
[\partial_+,\partial_-]=2\partial_0.$
\end{rem}

\begin{thm}\label{THM:tlmnpartialderivatives}
We have
\begin{eqnarray*}
  \partial_+ t^l_{mn} & = & -\sqrt{(l-n)(l+n+1)}\ t^l_{m,n+1}, \\
  \partial_- t^l_{mn} & = & -\sqrt{(l+n)(l-n+1)}\ t^l_{m,n-1}, \\
  \partial_0 t^l_{mn} & = & n\ t^l_{mn}, \\
  \Lap t^l_{mn} & = & -l(l+1)\ t^l_{mn}.
\end{eqnarray*}
\end{thm}

\begin{proof}
Formulae for $\partial_+,\partial_-,\partial_0$
follow from calculations in \cite[p. 141-142]{Vilenkin} .
Since
$$
  \Lap = - \partial_0^2 - (\partial_+\partial_-+\partial_-\partial_+)/2,
$$
we get
\begin{eqnarray*}
  \Lap t^l_{mn}
  & = & -n^2\ t^l_{mn} + \frac{1}{2} \left(
    \sqrt{(l+n)(l-n+1)}\ \partial_+ t^l_{m,n-1} \right. \\
    && \quad\quad\quad\quad\quad \left. +
    \sqrt{(l-n)(l+n+1)}\ \partial_- t^l_{m,n+1} \right) \\
  & = & \frac{-1}{2}
  \left( 2n^2 + \sqrt{(l+n)(l-n+1)} \sqrt{(l-(n-1))(l+(n-1)+1)}
    \right. \\
    && \quad\quad\quad \left.
    + \sqrt{(l-n)(l+n+1)} \sqrt{(l+(n+1))(l-(n+1)+1)} 
    \right) t^l_{mn} \\
  & = & \frac{-1}{2}
  \left( 2n^2 + (l+n)(l-n+1) + (l-n)(l+n+1) \right) t^l_{mn} \\
  & = & \frac{-1}{2} 
  \left(2n^2 + 2(l^2-n^2) + (l+n) + (l-n)\right) t^l_{mn} \\
  & = & -l(l+1)\ t^l_{mn}.
\end{eqnarray*}
\end{proof}

We can now calculate symbols of $\partial_+,\partial_-,\partial_0$
and of the Laplacian $\Lap$.
\begin{thm}\label{diffsymbols}
We have
\begin{eqnarray*}
  \sigma_{\partial_+}(x,l)_{mn} & = & -\sqrt{(l-n)(l+n+1)}\ 
  \delta_{m,n+1}
   =  -\sqrt{(l-m+1)(l+m)}\ \delta_{m-1,n},\\
  \sigma_{\partial_-}(x,l)_{mn} & = & -\sqrt{(l+n)(l-n+1)}\ 
  \delta_{m,n-1}
   =  -\sqrt{(l+m+1)(l-m)}\ \delta_{m+1,n},\\
  \sigma_{\partial_0}(x,l)_{mn} & = & n\ \delta_{mn}
   =  m\ \delta_{mn},\\
  \sigma_\Lap(x,l)_{mn} & = & -l(l+1)\ \delta_{mn},
\end{eqnarray*}
where $\delta_{mn}$ is the Kronecker delta:
$\delta_{mn}=1$ for $m=n$ and,
$\delta_{mn}=0$ otherwise.
\end{thm}

\begin{proof}
Let $e\in\SU2$ be the neutral element of $\SU2$ and let
$t^l$ be a unitary matrix representation of $\SU2$.
First we note that
$$
  \delta_{mn} = t^l(e)_{mn}
  = t^l(x^{-1}x)_{mn}
  = \sum_k t^l(x^{-1})_{mk}\ t^l(x)_{kn}
  = \sum_k \overline{t^l(x)_{km}}\ t^l(x)_{kn}.
$$
Similarly,
$\displaystyle
\delta_{mn} = \sum_k t^l(x)_{mk}\ \overline{t^l(x)_{nk}}$.
From this, formulae 
\eqref{EQ:su2-symbol-G-1}-\eqref{EQ:su2-symbol-G-2},
and Theorem~\ref{THM:tlmnpartialderivatives}
we get
\begin{eqnarray*}
   \sigma_{\partial_+}(x,l)_{mn}
   & = & \sum_k \overline{t^l_{km}(x)}\ \left(\partial_+ t^l_{kn}\right)(x)\\
   & = & -\sqrt{(l-n)(l+n+1)} \sum_k
   \overline{t^l_{km}(x)}\ t^l_{k,n+1}(x)\\
   & = & -\sqrt{(l-n)(l+n+1)}\ \delta_{m,n+1},
\end{eqnarray*}
and the case of $\sigma_{\partial_-}(x,l)$ is analogous.
Finally,
\begin{equation*}
   \sigma_{\partial_0}(x,l)_{mn}
    =  \sum_k \overline{t^l_{km}(x)}\ 
    \left(\partial_0 t^l_{kn}\right)(x)
    =  n \sum_k
   \overline{t^l_{km}(x)}\ t^l_{k,n}(x)
    =  n\ \delta_{m,n},
\end{equation*}
and similarly for $\Lap$, completing the proof.
\end{proof}

\begin{rem}\label{REM:su2-creation}
Notice that $\sigma_{\partial_0}(x,l)$ and $\sigma_\Lap(x,l)$
are diagonal matrices.
The non-zero elements reside
just above the diagonal of $\sigma_{\partial_+}(x,l)$, and
just below the diagonal of $\sigma_{\partial_-}(x,l)$.
Because of this operators $\partial_0$, $\partial_+$ and
$\partial_-$ may be called neutral, creation and
annihilation operators, respectively, and this explains
our preference to work with them rather than with $D_j$'s,
which have more non-zero entries.
\end{rem}

\section{Differences for symbols on $\SU2$}
\label{SEC:su2}

In this section we describe difference operators on $\SU2$
leading to symbol inequalities for symbols introduced in
Definition~\ref{DEF:su2-symbols-on-G}.
From Proposition~\ref{matrixelements} and 
Theorem~\ref{multiplication}
we recall the notation
\begin{eqnarray*}
  t^{1/2} & = &
  \begin{pmatrix}
    t_{--} & t_{-+} \\
    t_{+-} & t_{++}
    \end{pmatrix} 
   = 
  \begin{pmatrix}
    t^{1/2}_{-1/2,-1/2} & t^{1/2}_{-1/2,+1/2} \\
    t^{1/2}_{+1/2,-1/2} & t^{1/2}_{+1/2,+1/2}
    \end{pmatrix}.
\end{eqnarray*}
\begin{defn}\label{DEF:qs-su2}
For $q\in C^\infty(\SU2)$ and $f\in{\mathcal D}'(\SU2)$,
let $\triangle_q\widehat{f}(l) := \widehat{qf}(l)$.
We shall use abbreviations $\triangle_+=\triangle_{q_+}$,
$\triangle_-=\triangle_{q_-}$ and $\triangle_0=\triangle_{q_0}$,
where
\begin{eqnarray*}
  q_- &:=& t_{-+} = t^{1/2}_{-1/2,+1/2}, \\
  q_+ &:=& t_{+-} = t^{1/2}_{+1/2,-1/2}, \\
  q_0 &:=& t_{--}-t_{++} = 
  t^{1/2}_{-1/2,-1/2}-t^{1/2}_{+1/2,+1/2}.
\end{eqnarray*}
\end{defn}
Thus each trigonometric polynomial $q_+,q_-,q_0\in C^\infty(\SU2)$
vanishes at the neutral element $e\in \SU2$. In this sense  
trigonometric polynomials $q_- +q_+, q_- - q_+, q_0$ on $\SU2$
are analogues of polynomials $x_1, x_2,x_3$ in the
Euclidean space $\R^3$.

The aim now is to define difference operators acting on
symbols. For this purpose we may only look at symbols
independent of $x$ corresponding to right invariant operators
since the following construction is independent of $x$.
Thus, let $a=a(\xi)$ be a symbol as in
Definition~\ref{DEF:su2-symbols-on-G}. It follows that
$a=\widehat{s}$ for some right-convolution 
kernel $s\in {\mathcal D}'(\SU2)$ so that operator $Op(a)$
is given by
$$
  Op(a) f = f\ast s.
$$
Let us define ``difference operators''
$\triangle_+,\triangle_-,\triangle_0$
acting on symbol $a$ by
\begin{eqnarray}\label{differenceoperators1}
  \triangle_+ a &:=& \widehat{q_+\ s}, \\
  \triangle_- a &:=& \widehat{q_-\ s}, \\
  \triangle_0 a &:=& \widehat{q_0\ s}.
  \label{differenceoperators3}
\end{eqnarray}
We note that this construction is analogous to the one
producing usual derivatives in $\Rn$ or difference 
operators on the torus $\Tn$ (see \cite{RT08} for details).
On $\SU2$, to analyse the structure of these difference
operators, we first need to know how to multiply
functions $t^l_{mn}$ by $q_+,q_-,q_0$, and the
necessary
formulae are given in Theorem~\ref{multiplication}.

Let us now derive explicit expressions for the first
order difference operators $\triangle_+$, $\triangle_-$,
$\triangle_0$ defined in 
\eqref{differenceoperators1}-\eqref{differenceoperators3}. 
To abbreviate the notation, we will also
write $a^l_{nm}=a(x,l)_{nm}$, even if symbol $a(x,l)$
depends on $x$, keeping in mind that 
the following theorem holds pointwise in $x$. 

\begin{thm}\label{THM:differences}
The difference operators are given by
\begin{eqnarray*}
(\triangle_- a)^l_{nm} & = &
\frac{\sqrt{(l-m)(l+n)} }{2l+1}
  a^{l^-}_{n^- m^+} - 
  \frac{\sqrt{(l+m+1)(l-n+1)}}{2l+1} 
  a^{l^+}_{n^- m^+}, \\
 (\triangle_+ a)^l_{nm} & = &
\frac{\sqrt{(l+m)(l-n)}}{2l+1}
  a^{l^-}_{n^+ m^-} - 
  \frac{\sqrt{(l-m+1)(l+n+1)}}{2l+1} 
  a^{l^+}_{n^+ m^-}, \\
 (\triangle_0 a)^l_{nm} & = &
\frac{\sqrt{(l-m)(l-n)}}{2l+1}
  a^{l^-}_{n^+ m^+} + 
  \frac{\sqrt{(l+m+1)(l+n+1)}}{2l+1} 
  a^{l^+}_{n^+ m^+}- \\  
  & & -\frac{\sqrt{(l+m)(l+n)}}{2l+1}
  a^{l^-}_{n^- m^-} - 
  \frac{\sqrt{(l-m+1)(l-n+1)}}{2l+1} 
  a^{l^+}_{n^- m^-}, 
 \end{eqnarray*}
 where $k^\pm=k\pm\frac12,$ and satisfy commutator relations
\begin{equation}\label{EQ:sym-comm}
 [\triangle_0,\triangle_+]=[\triangle_0,\triangle_-]=
 [\triangle_-,\triangle_+]=0.
\end{equation} 
\end{thm}

\begin{proof}
Equalities \eqref{EQ:sym-comm} follow immediately
from \eqref{differenceoperators1}--\eqref{differenceoperators3}.
We can abbreviate $a(x,l)$ by $a(l)$ since none of the arguments
in the proof will act on the variable $x$.
Recall that by \eqref{EQ:su2-right-conv} we have 
$$
 a(x,l)_{nm} = a^l_{nm}
 = \widehat{s}(l)_{nm}
 = \int_\SU2 s(y)\ \overline{t^l_{mn}(y)}\ {\rm d}y,
$$
and
\begin{equation}\label{EQ:su2-kernel-symbol}
  s(x) = \sum_l (2l+1)\ {\rm Tr}\left( a(x,l)\ t^l(x) \right)
 = \sum_l (2l+1)\sum_{m,n}a^l_{nm}\ t^l_{mn}.
\end{equation} 
In the calculation below we will not worry about boundaries
of summations keeping in mind that we can always view
finite matrices as infinite ones simply by extending them
be zeros.
Recalling that $q_-=t_{-+}$ and using
Theorem~\ref{multiplication}, we can calculate
\begin{eqnarray*}
  q_-\ s
  & = & 
   \sum_l (2l+1)\sum_{m,n}a^l_{nm}\ q_-\ t^l_{mn} \\
  & = & \sum_l \sum_{m,n} a^l_{nm}
  \left[t^{l^+}_{m^- n^+} \sqrt{(l-m+1)(l+n+1)} 
     - t^{l^-}_{m^- n^+}\sqrt{(l+m)(l-n)} \right]  \\
  & = & \sum_l \sum_{m,n} t^l_{mn}
    \left[ a^{l^-}_{n^- m^+} \sqrt{(l-m)(l+n)} 
    - a^{l^+}_{n^- m^+} \sqrt{(l+m+1)(l-n+1)} \right].
\end{eqnarray*}
Since $\triangle_- a=\widehat{q_- s}$, we obtain
the desired formula for $\triangle_-$.
The calculation for $\triangle_+$ is analogous.
Finally, for $\triangle_0$, we calculate
\begin{eqnarray*}
  q_0\ s
  & = & 
   \sum_l (2l+1)\sum_{m,n}a^l_{nm}\ q_0\ t^l_{mn} \\
  & = & \sum_l \sum_{m,n} a^l_{nm} 
  \left[t^{l^+}_{m^- n^-} \sqrt{(l-m+1)(l-n+1)} 
  + t^{l^-}_{m^- n^-} \sqrt{(l+m)(l+n)} \right.  \\
  & & - \left.
  t^{l^+}_{m^+ n^+} \sqrt{(l+m+1)(l+n+1)} 
  - t^{l^-}_{m^+ n^+} \sqrt{(l-m)(l-n)} \right] \\
   & = & \sum_l \sum_{m,n} t^{l}_{mn}
   \left[ a^{l^-}_{n^+ m^+} \sqrt{(l-m)(l-n)}
  + a^{l^+}_{n^+ m^+} \sqrt{(l+m+1)(l+n+1)} \right.
  \\ & & \left. - a^{l^-}_{n^- m^-} \sqrt{(l+m)(l+n)} 
  - a^{l^+}_{n^- m^-} \sqrt{(l-m+1)(l-n+1)} \right].
\end{eqnarray*}
From this we obtain
the desired formula for $\triangle_0$ and the proof of 
Theorem~\ref{THM:differences} is complete.
\end{proof}

Let us now calculate higher order differences of
symbol $a\sigma_{\partial_0}$ which will be needed in the sequel.
\begin{thm}\label{THM:su2-Leibnitz}
For any $\alpha\in \Bbb N_0^3$, we have the formula
\begin{multline*}
\left[\triangle_+^{\alpha_1}\triangle_-^{\alpha_2}
\triangle_0^{\alpha_3} 
(a\sigma_{\partial_0})\right]^l_{nm} = \\
= (m-\alpha_1/2+\alpha_2/2)\left[
\triangle_+^{\alpha_1}\triangle_-^{\alpha_2} 
\triangle_0^{\alpha_3}
a\right]^l_{nm}
+\alpha_3\left[\overline{\triangle_0}
\triangle_+^{\alpha_1}\triangle_-^{\alpha_2} 
\triangle_0^{\alpha_3-1}
a\right]^l_{nm},
\end{multline*}
where $\overline{\triangle_0}$ is given by
\begin{eqnarray*}
(\overline{\triangle_0} a)^l_{nm} & = & 
  \frac12\left[
  \frac{\sqrt{(l-m)(l-n)}}{2l+1} 
   a^{l^-}_{n^+ m^+}+
  \frac{\sqrt{(l+m+1)(l+n+1)}}{2l+1} 
   a^{l^+}_{n^+ m^+}+ \right. \\
  &  & \left. +
 \frac{\sqrt{(l+m)(l+n)}}{2l+1} 
   a^{l^-}_{n^- m^-}+
  \frac{\sqrt{(l-m+1)(l-n+1)}}{2l+1} 
   a^{l^+}_{n^- m^-} 
  \right],
\end{eqnarray*}
and satisfies $[\triangle_0,\overline{\triangle_0}]=0$. 
\end{thm}

\begin{proof}
First we observe that we have
$$
 (a\ \sigma_{\partial_0})^l_{nm}
 =\sum_k a^l_{nk}\  k\ \delta_{km} = m\ a^l_{nm}.
$$
Then using Theorem~\ref{THM:differences}, we get
\begin{eqnarray*}
\triangle_-(a\sigma_{\partial_0})^l_{nm} & = &
 \frac{\sqrt{(l-m)(l+n)}}{2l+1} 
  m^+ a^{l^-}_{n^- m^+}-
  \frac{\sqrt{(l+m+1)(l-n+1)}}{2l+1} 
  m^+ a^{l^+}_{n^- m^+} \\
& = & (m^+ \triangle_- a)^l_{nm},
\end{eqnarray*}
and we can abbreviate this by writing
$\triangle_-(a \sigma_{\partial_0})=m^+\triangle_- a$.
Further, we have
\begin{eqnarray*} & &
\triangle_-(\triangle_-(a\sigma_{\partial_0}))^l_{nm} = \\ & = &
 \frac{\sqrt{(l-m)(l+n)}}{2l+1} 
  \left[\triangle_-(a\sigma_{\partial_0})\right]^{l^-}_{n^- m^+}-
  \frac{\sqrt{(l+m+1)(l-n+1)}}{2l+1}
 \left[\triangle_-(a\sigma_{\partial_0})\right]^{l^+}_{n^- m^+}
   \\ & = &
   \frac{\sqrt{(l-m)(l+n)}}{2l+1} 
  (m+1)(\triangle_- a)^{l^-}_{n^- m^+}-
  \frac{\sqrt{(l+m+1)(l-n+1)}}{2l+1}
  (m+1)(\triangle_- a)^{l^+}_{n^- m^+} \\
& = & (m+1) (\triangle_-^2 a)^l_{nm}.
\end{eqnarray*}
Continuing this calculation we can obtain
\begin{equation}\label{EQ:triangle-plus-k}
\left[\triangle_-^k(a\sigma_{\partial_0})\right]^l_{nm}=
(m+k/2)(\triangle_-^k a)^l_{nm}.
\end{equation} 
By Theorem~\ref{THM:differences} we also have
\begin{eqnarray*} & &
\left[\triangle_+(\triangle_-(a\sigma_{\partial_0}))
 \right]^l_{nm} = \\ & = &
 \left[\triangle_+\left(m^+ \triangle_- a\right)\right]^l_{nm}
 \\ & = & 
 \frac{\sqrt{(l+m)(l-n)}}{2l+1} 
  \left(m^+\triangle_- a\right)^{l^-}_{n^+ m^-}-
  \frac{\sqrt{(l-m+1)(l+n+1)}}{2l+1}
 \left(m^+\triangle_- a)\right)^{l^+}_{n^+ m^-}
   \\ 
& = & m (\triangle_+\triangle_- a)^l_{nm}.
\end{eqnarray*}
By induction we then get
\begin{equation}\label{EQ:triangle-plusmin-k}
\left[\triangle_+^{k_1}\triangle_-^{k_2}
(a\sigma_{\partial_0})\right]^l_{nm}=
(m-k_1/2+k_2/2)(\triangle_+^{k_1}\triangle_-^{k_2} a)^l_{nm}.
\end{equation} 
The situation with $\triangle_0$ is more complicated because
there are more terms. Using Theorem~\ref{THM:differences}
we have
\begin{eqnarray*} & & 
\triangle_0(a\sigma_{\partial_0})^l_{nm} = \\ & = &
 \frac{\sqrt{(l-m)(l-n)}}{2l+1} 
  ( m a)^{l^-}_{n^+ m^+}+
  \frac{\sqrt{(l+m+1)(l+n+1)}}{2l+1} 
  ( m a)^{l^+}_{n^+ m^+}- \\
  &  & -
 \frac{\sqrt{(l+m)(l+n)}}{2l+1} 
  ( m a)^{l^-}_{n^- m^-}-
  \frac{\sqrt{(l-m+1)(l-n+1)}}{2l+1} 
  ( m a)^{l^+}_{n^- m^-} \\
  & = & 
  \frac{\sqrt{(l-m)(l-n)}}{2l+1} 
  m^+ a^{l^-}_{n^+ m^+}+
  \frac{\sqrt{(l+m+1)(l+n+1)}}{2l+1} 
  m^+ a^{l^+}_{n^+ m^+}- \\
  &  & -
 \frac{\sqrt{(l+m)(l+n)}}{2l+1} 
  m^- a^{l^-}_{n^- m^-}-
  \frac{\sqrt{(l-m+1)(l-n+1)}}{2l+1} 
  m^- a^{l^+}_{n^- m^-} \\
  & = & m(\triangle_0 a)^l_{nm}+
  \frac12\left[
  \frac{\sqrt{(l-m)(l-n)}}{2l+1} 
   a^{l^-}_{n^+ m^+}+
  \frac{\sqrt{(l+m+1)(l+n+1)}}{2l+1} 
   a^{l^+}_{n^+ m^+}+ \right. \\
  &  & \left. +
 \frac{\sqrt{(l+m)(l+n)}}{2l+1} 
   a^{l^-}_{n^- m^-}+
  \frac{\sqrt{(l-m+1)(l-n+1)}}{2l+1} 
   a^{l^+}_{n^- m^-} 
  \right] \\
  & = & m(\triangle_0 a)^l_{nm} +
  (\overline{\triangle_0} a)^l_{nm},
\end{eqnarray*}
where $\overline{\triangle_0}$ is a weighted
averaging operator given by
\begin{eqnarray*}
(\overline{\triangle_0} a)^l_{nm} & = & 
  \frac12\left[
  \frac{\sqrt{(l-m)(l-n)}}{2l+1} 
   a^{l^-}_{n^+ m^+}+
  \frac{\sqrt{(l+m+1)(l+n+1)}}{2l+1} 
   a^{l^+}_{n^+ m^+}+ \right. \\
  &  & \left. +
 \frac{\sqrt{(l+m)(l+n)}}{2l+1} 
   a^{l^-}_{n^- m^-}+
  \frac{\sqrt{(l-m+1)(l-n+1)}}{2l+1} 
   a^{l^+}_{n^- m^-} 
  \right].
\end{eqnarray*}
We want to find a formula for $\triangle_0^k$, and for this
we first calculate
\begin{eqnarray*} & & 
\left[\triangle_0(\overline{\triangle_0} a)
 \right]^l_{nm} = \\ & = &
 \frac{\sqrt{(l-m)(l-n)}}{2l+1} 
  (\overline{\triangle_0} a)^{l^-}_{n^+ m^+}+
  \frac{\sqrt{(l+m+1)(l+n+1)}}{2l+1} 
  (\overline{\triangle_0} a)^{l^+}_{n^+ m^+}- \\
  &  & -
 \frac{\sqrt{(l+m)(l+n)}}{2l+1} 
  (\overline{\triangle_0} a)^{l^-}_{n^- m^-}-
  \frac{\sqrt{(l-m+1)(l-n+1)}}{2l+1} 
  (\overline{\triangle_0} a)^{l^+}_{n^- m^-} \\
  & = & 
  \frac{\sqrt{(l-m)(l-n)}}{2l+1} \frac12 \left[
  \frac{\sqrt{(l^--m^+)(l^--n^+)}}{2l^-+1} 
   a^{l^{--}}_{n^{++} m^{++}}+ \right. \\ & & +
  \frac{\sqrt{(l^-+m^++1)(l^-+n^++1)}}{2l^-+1} 
  a^{l^{-+}}_{n^{++} m^{++}}+  \\
  &  & + \left.
  \frac{\sqrt{(l^-+m^+)(l^-+n^+)}}{2l^-+1} 
   a^{l^{--}}_{n^{+-} m^{+-}}+
  \frac{\sqrt{(l^--m^++1)(l^--n^++1)}}{2l^-+1} 
  a^{l^{-+}}_{n^{+-} m^{+-}} 
  \right]+ \\
   & & + 
  \frac{\sqrt{(l+m+1)(l+n+1)}}{2l+1} \frac12 \frac{1}{2l^++1}
  \left[
  \sqrt{(l^+-m^+)(l^+-n^+)} 
   a^{l^{+-}}_{n^{++} m^{++}}+ \right. \\ & & +
  \sqrt{(l^++m^++1)(l^++n^++1)} 
  a^{l^{++}}_{n^{++} m^{++}}+  \\
  &  & + \left.
  \sqrt{(l^++m^+)(l^++n^+)} 
   a^{l^{+-}}_{n^{+-} m^{+-}}+
  \sqrt{(l^+-m^++1)(l^+-n^++1)} 
  a^{l^{++}}_{n^{+-} m^{+-}} 
  \right]- 
  \\
   & & -
  \frac{\sqrt{(l+m)(l+n)}}{2l+1} \frac12 \frac{1}{2l^-+1}
  \left[
  \sqrt{(l^--m^-)(l^--n^-)} 
   a^{l^{--}}_{n^{-+} m^{-+}}+ \right. \\ & & +
  \sqrt{(l^-+m^-+1)(l^-+n^-+1)} 
  a^{l^{-+}}_{n^{-+} m^{-+}}+  \\
  &  & + \left.
  \sqrt{(l^-+m^-)(l^-+n^-)} 
   a^{l^{--}}_{n^{--} m^{--}}+
  \sqrt{(l^--m^-+1)(l^--n^-+1)} 
  a^{l^{-+}}_{n^{--} m^{--}} 
  \right]- 
  \\
   & & -
  \frac{\sqrt{(l-m+1)(l-n+1)}}{2l+1} \frac12 \frac{1}{2l^++1}
  \left[
  \sqrt{(l^+-m^-)(l^+-n^-)} 
   a^{l^{+-}}_{n^{-+} m^{-+}}+ \right. \\ & & +
  \sqrt{(l^++m^-+1)(l^++n^-+1)} 
  a^{l^{++}}_{n^{-+} m^{-+}}+  \\
  &  & + \left.
  \sqrt{(l^++m^-)(l^++n^-)} 
   a^{l^{+-}}_{n^{--} m^{--}}+
  \sqrt{(l^+-m^-+1)(l^+-n^-+1)} 
  a^{l^{++}}_{n^{--} m^{--}} 
  \right].
\end{eqnarray*}
From this we get
\begin{eqnarray*} & &
\left[\triangle_0(\overline{\triangle_0} a)
 \right]^l_{nm} = \\ & = &
  \frac{\sqrt{(l-m)(l-n)}}{2l+1} \frac12 \frac{1}{2l} \left[
  \sqrt{(l-m-1)(l-n-1)} 
   a^{l^{--}}_{n^{++} m^{++}}+ \right. \\ & & +
  \sqrt{(l+m+1)(l+n+1)} 
  a^l_{n^{++} m^{++}}+  \\
  &  & + \left.
  \sqrt{(l+m)(l+n)} 
   a^{l^{--}}_{nm}+
  \sqrt{(l-m)(l-n)} 
  a^l_{nm} 
  \right]+ \\
   & & + 
  \frac{\sqrt{(l+m+1)(l+n+1)}}{2l+1} \frac12 \frac{1}{2l+2}
  \left[
  \sqrt{(l-m)(l-n)} 
   a^l_{n^{++} m^{++}}+ \right. \\ & & +
  \sqrt{(l+m+2)(l+n+2)} 
  a^{l^{++}}_{n^{++} m^{++}}+  \\
  &  & + \left.
  \sqrt{(l+m+1)(l+n+1)} 
   a^l_{nm}+
  \sqrt{(l-m+1)(l-n+1)} 
  a^{l^{++}}_{nm} 
  \right]- 
  \\
   & & -
  \frac{\sqrt{(l+m)(l+n)}}{2l+1} \frac12 \frac{1}{2l}
  \left[
  \sqrt{(l-m)(l-n)} 
   a^{l^{--}}_{nm}+ \right. \\ & & +
  \sqrt{(l+m)(l+n)} 
  a^l_{nm}+  \\
  &  & + \left.
  \sqrt{(l+m-1)(l+n-1)} 
   a^{l^{--}}_{n^{--} m^{--}}+
  \sqrt{(l-m+1)(l-n+1)} 
  a^l_{n^{--} m^{--}} 
  \right]- 
  \\
   & & -
  \frac{\sqrt{(l-m+1)(l-n+1)}}{2l+1} \frac12 \frac{1}{2l+2}
  \left[
  \sqrt{(l-m+1)(l-n+1)} 
   a^l_{nm}+ \right. \\ & & +
  \sqrt{(l+m+1)(l+n+1)} 
  a^{l^{++}}_{nm}+  \\
  &  & + \left.
  \sqrt{(l+m)(l+n)} 
   a^l_{n^{--} m^{--}}+
  \sqrt{(l-m+2)(l-n+2)} 
  a^{l^{++}}_{n^{--} m^{--}} 
  \right],
  \end{eqnarray*}
and we can note that here pairs of terms 
with $a^{l^{--}}_{nm}$, $a^{l^{++}}_{nm}$
cancel, and also four terms with
$a^l_{nm}$ cancel in view of the identity
\begin{eqnarray*}
\frac{(l-m)(l-n)}{(2l+1)(2l)}+
\frac{(l+m+1)(l+n+1)}{(2l+1)(2l+2)}-
\frac{(l+m)(l+n)}{(2l+1)(2l)}-
\frac{(l-m+1)(l-n+1)}{(2l+1)(2l+2)} & & \\
= \frac{-2l(m+n)}{(2l+1)(2l)}+
\frac{(2l+2)(m+n)}{(2l+1)(2l+2)} =0.  & &
\end{eqnarray*}
Calculating in the other direction, we get
\begin{eqnarray*} & &
\left[\overline{\triangle_0}({\triangle_0} a)
 \right]^l_{nm} = \\ & = &
 \frac12 \frac{\sqrt{(l-m)(l-n)}}{2l+1} 
  ({\triangle_0} a)^{l^-}_{n^+ m^+}+
  \frac12 \frac{\sqrt{(l+m+1)(l+n+1)}}{2l+1} 
  ({\triangle_0} a)^{l^+}_{n^+ m^+}+ \\
  &  & +
 \frac12 \frac{\sqrt{(l+m)(l+n)}}{2l+1} 
  ({\triangle_0} a)^{l^-}_{n^- m^-}
  +\frac12 \frac{\sqrt{(l-m+1)(l-n+1)}}{2l+1} 
  ({\triangle_0} a)^{l^+}_{n^- m^-} \\
  & = & 
  \frac{\sqrt{(l-m)(l-n)}}{2l+1} \frac12 \frac{1}{2l^-+1}\left[
  \sqrt{(l^--m^+)(l^--n^+)} 
   a^{l^{--}}_{n^{++} m^{++}}+ \right. \\ & & +
  \sqrt{(l^-+m^++1)(l^-+n^++1)} 
  a^{l^{-+}}_{n^{++} m^{++}}-  \\
  &  & - \left.
  \sqrt{(l^-+m^+)(l^-+n^+)} 
   a^{l^{--}}_{n^{+-} m^{+-}}-
  \sqrt{(l^--m^++1)(l^--n^++1)} 
  a^{l^{-+}}_{n^{+-} m^{+-}} 
  \right]+ \\
   & & +
  \frac{\sqrt{(l+m+1)(l+n+1)}}{2l+1} \frac12 \frac{1}{2l^++1}
  \left[
  \sqrt{(l^+-m^+)(l^+-n^+)} 
   a^{l^{+-}}_{n^{++} m^{++}}+ \right. \\ & & +
  \sqrt{(l^++m^++1)(l^++n^++1)} 
  a^{l^{++}}_{n^{++} m^{++}}-  \\
  &  & - \left.
  \sqrt{(l^++m^+)(l^++n^+)} 
   a^{l^{+-}}_{n^{+-} m^{+-}}-
  \sqrt{(l^+-m^++1)(l^+-n^++1)} 
  a^{l^{++}}_{n^{+-} m^{+-}} 
  \right]+ 
  \\
   & & +
  \frac{\sqrt{(l+m)(l+n)}}{2l+1} \frac12 \frac{1}{2l^-+1}
  \left[
  \sqrt{(l^--m^-)(l^--n^-)} 
   a^{l^{--}}_{n^{-+} m^{-+}}+ \right. \\ & & +
  \sqrt{(l^-+m^-+1)(l^-+n^-+1)} 
  a^{l^{-+}}_{n^{-+} m^{-+}}-  \\
  &  & - \left.
  \sqrt{(l^-+m^-)(l^-+n^-)} 
   a^{l^{--}}_{n^{--} m^{--}}-
  \sqrt{(l^--m^-+1)(l^--n^-+1)} 
  a^{l^{-+}}_{n^{--} m^{--}} 
  \right]+ 
  \\
   & & +
  \frac{\sqrt{(l-m+1)(l-n+1)}}{2l+1} \frac12 \frac{1}{2l^++1}
  \left[
  \sqrt{(l^+-m^-)(l^+-n^-)} 
   a^{l^{+-}}_{n^{-+} m^{-+}}+ \right. \\ & & +
  \sqrt{(l^++m^-+1)(l^++n^-+1)} 
  a^{l^{++}}_{n^{-+} m^{-+}}-  \\
  &  & - \left.
  \sqrt{(l^++m^-)(l^++n^-)} 
   a^{l^{+-}}_{n^{--} m^{--}}-
  \sqrt{(l^+-m^-+1)(l^+-n^-+1)} 
  a^{l^{++}}_{n^{--} m^{--}} 
  \right].
\end{eqnarray*}
From this we get
\begin{eqnarray*} & &
\left[\overline{\triangle_0}({\triangle_0} a)
 \right]^l_{nm} = \\ & = &
  \frac{\sqrt{(l-m)(l-n)}}{2l+1} \frac12 \frac{1}{2l} \left[
  \sqrt{(l-m-1)(l-n-1)} 
   a^{l^{--}}_{n^{++} m^{++}}+ \right. \\ & & +
  \sqrt{(l+m+1)(l+n+1)} 
  a^l_{n^{++} m^{++}}-  \\
  &  & - \left.
  \sqrt{(l+m)(l+n)} 
   a^{l^{--}}_{nm}-
  \sqrt{(l-m)(l-n)} 
  a^l_{nm} 
  \right]+ \\
   & & + 
  \frac{\sqrt{(l+m+1)(l+n+1)}}{2l+1} \frac12 \frac{1}{2l+2}
  \left[
  \sqrt{(l-m)(l-n)} 
   a^l_{n^{++} m^{++}}+ \right. \\ & & +
  \sqrt{(l+m+2)(l+n+2)} 
  a^{l^{++}}_{n^{++} m^{++}}-  \\
  &  & - \left.
  \sqrt{(l+m+1)(l+n+1)} 
   a^l_{nm}-
  \sqrt{(l-m+1)(l-n+1)} 
  a^{l^{++}}_{nm} 
  \right]+
  \\
   & & +
  \frac{\sqrt{(l+m)(l+n)}}{2l+1} \frac12 \frac{1}{2l}
  \left[
  \sqrt{(l-m)(l-n)} 
   a^{l^{--}}_{nm}+ \right. \\ & & +
  \sqrt{(l+m)(l+n)} 
  a^l_{nm}-  \\
  &  & - \left.
  \sqrt{(l+m-1)(l+n-1)} 
   a^{l^{--}}_{n^{--} m^{--}}-
  \sqrt{(l-m+1)(l-n+1)} 
  a^l_{n^{--} m^{--}} 
  \right]+
  \\
   & & +
  \frac{\sqrt{(l-m+1)(l-n+1)}}{2l+1} \frac12 \frac{1}{2l+2}
  \left[
  \sqrt{(l-m+1)(l-n+1)} 
   a^l_{nm}+ \right. \\ & & +
  \sqrt{(l+m+1)(l+n+1)} 
  a^{l^{++}}_{nm}-  \\
  &  & - \left.
  \sqrt{(l+m)(l+n)} 
   a^l_{n^{--} m^{--}}-
  \sqrt{(l-m+2)(l-n+2)} 
  a^{l^{++}}_{n^{--} m^{--}} 
  \right] 
\end{eqnarray*}
and we can note that here terms $a^l_{nm}$,
$a^{l^{--}}_{nm}$ and $a^{l^{++}}_{nm}$ cancel again.
From these calculations we obtain
$$\overline{\triangle_0} {\triangle_0} a=
{\triangle_0}\overline{\triangle_0} a.$$
Then we can easily see that
$$
 \triangle_0^2(ma)=\triangle_0(m\triangle_0 a+
 \overline{\triangle_0} a)=
 m\triangle_0^2 a+2 \overline{\triangle_0} \triangle_0 a,
$$
and, moreover,
$$
 \triangle_0^k(ma)=
 m\triangle_0^k a+k \overline{\triangle_0} \triangle_0^{k-1} a.
$$
Let us now apply this to \eqref{EQ:triangle-plusmin-k}.
Using commutativity of $\triangle_0, \triangle_+$ and
$\triangle_-$ from Theorem~\ref{THM:differences}, we get
\begin{eqnarray*} & &
\left[\triangle_+^{k_1}\triangle_-^{k_2}\triangle_0^{k_3} 
(a\sigma_{\partial_0})\right]^l_{nm} = \\ & = &
\left[\triangle_0^{k_3}\triangle_+^{k_1}\triangle_-^{k_2} 
(a\sigma_{\partial_0})\right]^l_{nm} \\
& = & 
\left[\triangle_0^{k_3}\left((m-k_1/2+k_2/2)
\triangle_+^{k_1}\triangle_-^{k_2} a\right)\right]^l_{nm}
\\ & = & 
\left[\triangle_0^{k_3}\left(m
\triangle_+^{k_1}\triangle_-^{k_2} a\right)\right]^l_{nm}
-\left[\triangle_0^{k_3}\left((k_1/2-k_2/2)
\triangle_+^{k_1}\triangle_-^{k_2} a\right)\right]^l_{nm}
\\ & = & 
m\left[\triangle_0^{k_3}
\triangle_+^{k_1}\triangle_-^{k_2} a\right]^l_{nm}
+k_3\left[\overline{\triangle_0}\triangle_0^{k_3-1}
\triangle_+^{k_1}\triangle_-^{k_2} a\right]^l_{nm}-
\\ & & -(k_1/2-k_2/2)\left[\triangle_0^{k_3}
\triangle_+^{k_1}\triangle_-^{k_2} a\right]^l_{nm}
\\ & = & 
(m-k_1/2+k_2/2)\left[
\triangle_+^{k_1}\triangle_-^{k_2} \triangle_0^{k_3}
a\right]^l_{nm}
+k_3\left[\overline{\triangle_0}
\triangle_+^{k_1}\triangle_-^{k_2} \triangle_0^{k_3-1}
a\right]^l_{nm},
\end{eqnarray*}
completing the proof.
\end{proof}

We now collect some properties of first-order differences.
\begin{thm}\label{THM:su2-first-order-symb}
We have 
\begin{equation}\label{EQ:su2-diffs1}
    \sigma_I = 
  \triangle_+ \sigma_{\partial_+} =
  \triangle_- \sigma_{\partial_-} =
  \triangle_0 \sigma_{\partial_0}.
\end{equation} 
If $\mu,\nu\in\{+,-,0\}$ are such that $\mu\not=\nu$, then
\begin{equation}\label{EQ:su2-diffs2}
  \triangle_\mu \sigma_{\partial_\nu} = 0,
\end{equation}
and for every $\nu\in\{+,-,0\}$, we have
\begin{equation}\label{EQ:su2-diffs3}
  \triangle_\nu \sigma_I(x) = 0.
\end{equation}
Moreover, if $\Lap$ is the bi-invariant Laplacian, then
\begin{equation}\label{EQ:su2-diffs4}
  \triangle_+ \sigma_{\Lap} = -\sigma_{\partial_-},\quad
  \triangle_- \sigma_{\Lap} = -\sigma_{\partial_+},\quad
  \triangle_0 \sigma_{\Lap} = -2\sigma_{\partial_0}.
\end{equation}
\end{thm}

\begin{proof}
Let us prove \eqref{EQ:su2-diffs1}.
From Theorem~\ref{multiplication}
we get an expression for $q_+ t^l_{mn}=t_{+-} t^l_{mn}$,
which is used in the following calculation
together with \eqref{EQ:su2-kernel-symbol} and
Theorem~\ref{diffsymbols}:
\begin{eqnarray*}
  q_+\ s_{\partial_+}
  & = & q_+ \sum_l (2l+1) \sum_{m,n} 
  \sigma_{\partial_+}(l)_{mn}\ t^l_{nm} \\
  & = & \sum_l \sum_{n} \sigma_{\partial_+}(l)_{n+1,n}
  \ (2l+1)\ q_+\ t^l_{n,n+1} \\
  & = & \sum_l \sum_n - \left(\sqrt{(l-n)(l+n+1)}\right)^2
  \left( t^{l^+}_{n^+ n^+} - t^{l^-}_{n^+ n^+} \right) \\
  & = & \sum_l (2l+1) \sum_k t^l_{k k} \\
  & = & s_I.
\end{eqnarray*}
Hence
$\triangle_+ \sigma_{\partial_+} = \sigma_I$.
Similarly, we can show that 
$\triangle_- \sigma_{\partial_-} = \sigma_I$
and that $\triangle_0 \sigma_{\partial_0} = \sigma_I$.
Let us now prove \eqref{EQ:su2-diffs2}. We have
\begin{eqnarray*}
  q_+\ s_{\partial_-}
  & = & q_+ \sum_l (2l+1) \sum_{m,n}
  \sigma_{\partial_-}(l)_{mn}\ t^l_{nm} \\
  & = & \sum_l \sum_{n}
  \sigma_{\partial_-}(l)_{n-1,n}\ (2l+1)\ q_+\ t^l_{n,n-1} \\
  & = & \sum_l \sum_n - \sqrt{(l+n)(l-n+1)} \\
  && \hskip 10mm
  \left( \sqrt{(l+n+1)(l-n+2)}\ t^{l+1/2}_{n+1/2,n-3/2} \right. \\
   && \left. \hskip 10mm
     - \sqrt{(l-n)(l+n-1)}\ t^{l-1/2}_{n+1/2,n-3/2} \right) \\
  & = & \sum_l \sum_n t^{l+1/2}_{n+1/2,n-3/2} \\
  && \hskip 10mm \left(
    -\sqrt{(l+n+1)(l-n+2)} \sqrt{(l-n+1)(l+n)} \right. \\
  && \left. \hskip 10mm
    + \sqrt{(l+n)(l-n+1)} \sqrt{(l+n+1)(l-n+2)} \right) \\    
  & = & 0.
\end{eqnarray*}
Analogously, one can readily show the rest of
\eqref{EQ:su2-diffs2}.
Let us now prove \eqref{EQ:su2-diffs3}. We have
\begin{eqnarray*}
  q_-\ s_I & = & q_- \sum_l (2l+1) \sum_{m,n} t^l_{mn} \\
  & = & \sum_l \sum_n (2l+1)\ q_-\ t^l_{nn} \\
  & = & \sum_l \sum_n \left(
    \sqrt{(l-n+1)(l+n+1)}\ t^{l^+}_{n^-,n^+}
    -\sqrt{(l+n)(l-n)}\ t^{l^-}_{n^-,n^+} \right) \\
  & = & \sum_l \sum_n t^{l^-}_{n^-,n^+} \left(
    \sqrt{(l-n)(l+n)}
    -\sqrt{(l+n)(l-n)} \right) \\
  & = & 0.
\end{eqnarray*}
Analogously, we have
$q_+\ s_I  = q_0\ s_I = 0$
which proves proves \eqref{EQ:su2-diffs3}.
Let us finally prove \eqref{EQ:su2-diffs4}. 
Since
$$
  \sigma_{\Lap}(x,l)_{mn} = -l(l+1)\ \delta_{mn}
$$
by Theorem~\ref{diffsymbols},
we get
\begin{eqnarray*}
  q_-\ s_{-\Lap}
  & = & q_-\ \sum_l (2l+1)\sum_{m,n}
  \sigma_{-\Lap}(x,l)_{mn}\ t^l_{nm} \\
  & = & \sum_l (2l+1) \sum_n l(l+1)\ q_-\ t^l_{nn}(y) \\
  & = & \sum_l \sum_n l(l+1) \left(
    +\sqrt{(l-n+1)(l+n+1)}\ t^{l^+}_{n^-,n^+} \right. \\
  && \left. \hskip 25mm
    -\sqrt{(l+n)(l-n)}\ t^{l^-}_{n^-,n^+} \right) \\
  & = & \sum_l \sum_n t^{l^+}_{n^-,n^+} \left(
    + l(l+1) \sqrt{(l-n+1)(l+n+1)} \right. \\
  && \left. \hskip 25mm
    - (l+1)(l+2) \sqrt{(l+n+1)(l-n+1)} \right) \\
  & = & \sum_l \sum_n -2(l+1) \sqrt{(l+n+1)(l-n+1)}\ 
  t^{l^+}_{n^-,n^+} \\
  & = & \sum_l (2l+1) \sum_n -\sqrt{(l+n)(l-n+1)}\ 
  t^{l}_{n-1,n} \\
  & = & s_{\partial_+}.
\end{eqnarray*}
Analogously, one can readily show that
$ q_+\ s_{-\Lap}=s_{\partial_-}$ and that
$ q_0\ s_{-\Lap}= 2 s_{\partial_0}$, completing
the proof.
\end{proof}

\begin{rem}\label{REM:differencesondiffsymbols}
In Theorem~\ref{THM:su2-first-order-symb}
we applied the differences on the symbols
of specific differential operators on $\SU2$.
In general, on a compact Lie group $G$,
a difference operator of order $|\gamma|$
applied to a symbol of a partial differential operator of order $N$
gives a symbol of order $N-|\gamma|$.
More precisely, let
\begin{equation*}\label{diffopD}
  D = \sum_{|\alpha|\leq N} c_\alpha(x)\ \partial_x^\alpha
\end{equation*}
be a partial differential operator
with coefficients $c_\alpha\in C^\infty(G)$.
For $q\in C^\infty(G)$ such that $q(e)=0$, we define 
difference operator
$\triangle_q$ acting on symbols by
\begin{equation*}
  \triangle_q \widehat{f}(\xi)
  := \widehat{qf}(\xi).
\end{equation*}
Then we obtain
\begin{equation*}
  \triangle_q \sigma_D(x,\xi)
  = \sum_{|\alpha|\leq N} c_\alpha(x)
  \sum_{\beta\leq\alpha} {\alpha\choose\beta} (-1)^{|\beta|}
  \ (\partial_x^\beta q)(e)\ \sigma_{\partial_x^{\alpha-\beta}}(x,\xi),
\end{equation*}
which is a symbol of a partial differential operator
of order at most $N-1$.
\end{rem}

\section{Taylor expansion on Lie groups}
\label{SEC:Taylor-expansion}

As Taylor polynomial expansions are useful in obtaining 
symbolic calculus
on $\Bbb R^n$, we would like to have analogous expansions on 
group $G$. Here, Taylor expansion formula on $G$ will be
obtained by embedding $G$ into some $\R^m$, using the
Taylor expansion formula in $\R^m$, and then restricting
it back to $G$.

Let $U\subset\Bbb R^m$ be an open neighbourhood of some point 
$\vec{e}\in\Bbb R^m$.
The $N$th order Taylor polynomial $P_N f:\Bbb R^m\to\Bbb C$
of $f\in C^\infty(U)$ at $\vec e$ is given by
\begin{equation*}
  P_N f(\vec x) = \sum_{\alpha\in\Bbb N_0^m:\ |\alpha|\leq N}
  \frac{1}{\alpha!}\ (\vec x-\vec e)^\alpha\ 
  \partial_x^\alpha f(\vec e).
\end{equation*}
Then the remainder $E_N f := f-P_N f$ satisfies
\begin{equation*}
  E_N f(\vec x) = \sum_{|\alpha|=N+1} 
  (\vec x-\vec e)^\alpha\ 
  f_\alpha(\vec x)
\end{equation*}
for some functions $f_\alpha\in C^\infty(U)$.
In particular,
$$
  E_N f(\vec x) = {\mathcal O}(\|\vec x-\vec e\|^{N+1})\quad 
  {\rm as}\quad \vec x\to \vec e.
$$
Let $G$ be a compact Lie group;
we would like to approximate a smooth function $u:G\to\Bbb C$
using a Taylor polynomial type expansion
nearby the neutral element $e\in G$.
We may assume that $G$ is a closed subgroup
of ${\rm GL}(n,\Bbb R)\subset\Bbb R^{n\times n}$,
the group of real invertible $(n\times n)$-matrices,
and thus a closed submanifold of the Euclidean space
of dimension $m=n^2$. This embedding of $G$ into $\R^m$
will be denoted by $x\mapsto \vec x$, and the image of
$G$ under this embedding will be still denoted by $G$.
Also, if $x\in G$, we may still write $x$ for 
$\vec x$ to simplify the
notation.
Let $U\subset\Bbb R^m$ be a small enough
open neighbourhood of $G\subset\Bbb R^m$
such that for each $\vec x\in U$ there exists
a unique nearest point $p(\vec x)\in G$ 
(with respect to the Euclidean distance).
For $u\in C^\infty(G)$ define $f\in C^\infty(U)$ by
$$
  f(\vec x) := u(p(\vec x)).
$$
The effect is that $f$ is constant in the directions
perpendicular to $G$.
As above, we may define the Euclidean Taylor polynomial
$P_N f:\Bbb R^m\to\Bbb C$ at $e\in G\subset\Bbb R^m$.
Let us define $P_N u:G\to\Bbb C$ as the restriction,
$$
  P_N u:=P_N f|_G.
$$
We call $P_N u\in C^\infty(G)$
a {\it Taylor polynomial of $u$ of order $N$ at $e\in G$}.
Then for $x\in G$, we have
$$
  u(x)-P_N u(x) = \sum_{|\alpha|=N+1} u_\alpha(x)\ (x-e)^\alpha
$$
for some functions $u_\alpha\in C^\infty(G)$, where
we set $(x-e)^\alpha:=(\vec x-\vec e)^\alpha$.
Taylor polynomials on $G$ are given by 
\begin{equation*}\label{EQ:Taylor-pol-G}
  P_N u(x) = \sum_{|\alpha|\leq N}
  \frac{1}{\alpha!}\ (x-e)^\alpha\ 
  \partial_x^{(\alpha)} u(e),
\end{equation*}
where we set $\partial_x^{(\alpha)} u(e)
:=\partial_x^{\alpha} f(\vec e)$.

Let us now consider especially $G={\rm SU}(2)$.
Recall the quaternionic identification
$$
  \left( x_0{\bf 1}+x_1{\bf i}+x_2{\bf j}+x_3{\bf k}
  \mapsto (x_0,x_1,x_2,x_3) \right):\Bbb H\to\Bbb R^4.
$$
Moreover, there is the identification 
$\Bbb H\supset\Bbb S^3\cong{\rm SU}(2)$,
$$
  \vec x=(x_0,x_1,x_2,x_3) \mapsto
  \begin{pmatrix} x_0 + {\rm i}x_3 & x_1 + {\rm i}x_2 \\
    -x_1 + {\rm i}x_2 & x_0 - {\rm i}x_3 \end{pmatrix}
  = \begin{pmatrix} x_{11} & x_{12} \\ x_{21} & x_{22} \end{pmatrix}
  =x.
$$
Hence we identify $(1,0,0,0)\in\Bbb R^4$
with the neutral element of ${\rm SU}(2)$.
Notice that
\begin{eqnarray*}
  q_+(x) & = & x_{12}\ =\ x_1+{\rm i}x_2, \\
  q_-(x) & = & x_{21}\ =\ -x_1+{\rm i}x_2, \\
  q_0(x) & = & x_{11}-x_{22}\ =\ 2{\rm i}x_3.
\end{eqnarray*}
A function $u\in C^\infty(\Bbb S^3)$
can be extended to $f\in C^\infty(U)=C^\infty(\Bbb R^4\setminus\{0\})$
by
$$
  f(\vec x) := u(\vec x/\|\vec x\|).
$$
Then we obtain $P_N u\in C^\infty(\Bbb S^3)$,
$$
  P_N u(\vec x) := \sum_{|\alpha|\leq N} \frac{1}{\alpha!}
  \ \left(\vec x-\vec e\right)^{\alpha}\ 
  \partial_x^\alpha f(\vec e),
$$
where $\vec e=(1,0,0,0)$. Expressing this in terms of $x\in\SU2$,
we obtain Taylor polynomials for $x\in\SU2$:
$$
  P_N u(x) = \sum_{|\alpha|\leq N} \frac{1}{\alpha!}
  \ \left(x-e\right)^{\alpha}\ 
  \partial_x^{(\alpha)} u(e),
$$
where we write $\partial_x^{(\alpha)} u(e)
=\partial_x^{\alpha} f(\vec e)$, and where
\begin{multline*}
(x-e)^\alpha=(\vec x-\vec e)^\alpha=
(x_0-1)^{\alpha_1} x_1^{\alpha_2}x_2^{\alpha_3} 
x_3^{\alpha_4}= \\
= \p{\frac{x_{11}+x_{22}}{2}-1}^{\alpha_1}
\p{\frac{x_{12}-x_{21}}{2}}^{\alpha_2}
\p{\frac{x_{12}+x_{21}}{2\irm}}^{\alpha_3}
\p{\frac{x_{11}-x_{22}}{2\irm}}^{\alpha_4}.
\end{multline*}

\section{Properties of global pseudo-differential symbols}
\label{SEC:psdosymbols}

In this section, we study the global symbols
of pseudo-differential operators on compact Lie groups.
We also derive elements of the calculus
in more general classes of symbols, and prove 
the Sobolev boundedness Theorem~\ref{THM:su2-Sobolev}.

As explained in Section~\ref{SEC:Taylor-expansion}, 
smooth functions on a group $G$ can be approximated
by Taylor polynomial type expansions.
More precisely,
there exist partial differential operators $\partial_x^{(\alpha)}$
of order $|\alpha|$ on $G$ such that for every $u\in C^\infty(G)$
we have
\begin{equation}\label{EQ:Taylor-exp}
 u(x) =  \sum_{|\alpha|\leq N} \frac{1}{\alpha!}
  \ q_\alpha(x^{-1})\ \partial_x^{(\alpha)} u(e)
  + \sum_{|\alpha|=N+1}  q_\alpha(x^{-1})\ u_\alpha(x)
 \sim  \sum_{\alpha\geq 0} \frac{1}{\alpha!}
  \ q_\alpha(x^{-1})\ \partial_x^{(\alpha)} u(e)
\end{equation} 
in a neighbourhood of $e\in G$,
where $u_\alpha\in C^\infty(G)$,
and $q_\alpha\in C^\infty(G)$ satisfy
$q_{\alpha+\beta} = q_\alpha q_\beta$.
Moreover, here $q_0\equiv 1$, and $q_\alpha(e)=0$ if 
$|\alpha|\geq 1$.
Let us define difference operators $\triangle_\xi^\alpha$
acting on Fourier coefficients by
$\triangle_\xi^\alpha\widehat{f}(\xi) := \widehat{q_\alpha f}(\xi)$.
Notice that
$\triangle_\xi^{\alpha+\beta}=\triangle_\xi^\alpha\triangle_\xi^\beta$.

\begin{rem}
The technical
choice of writing $q_\alpha(x^{-1})$ in \eqref{EQ:Taylor-exp}
is dictated by our desire to make asymptotic formulae
in Theorems \ref{THM:composition} and \ref{THM:adjoint}
look similar to the familiar Euclidean formulae, and by
an obvious freedom in selecting among different forms of
Taylor polynomials $q_\alpha$. 
For example, on $\SU2$, if we work with operators
$\Delta_+, \Delta_-, \Delta_0$ defined in 
\eqref{differenceoperators1}-\eqref{differenceoperators3},
we can choose the form of the Taylor expansion 
\eqref{EQ:Taylor-exp} adapted to functions $q_+,q_-,q_0$.
Here we can observe that $q_+(x^{-1})=-q_-(x),$
$q_-(x^{-1})=-q_+(x),$ $q_0(x^{-1})=-q_0(x),$ so that
for $|\alpha|=1$ functions $q_\alpha(x)$ and $q_\alpha(x^{-1})$
are linear combinations of $q_+, q_-, q_0$.
\end{rem}

Let $\{Y_j\}_{j=1}^{{\rm dim}(G)}$
be a basis for the Lie algebra of $G$,
and let $\partial_j$ be the left-invariant vector fields
corresponding to $Y_j$.
For $\beta\in\Bbb N_0^n$, let us denote
$\partial^\beta = \partial_1^{\beta_1}\cdots\partial_n^{\beta_n}$.

For a compact closed manifold $M$,
let ${\mathcal A}_0^m(M)$ denote the set of 
those continuous linear operators
$A:C^\infty(M)\to C^\infty(M)$ which are bounded 
from $H^m(M)$ to $L^2(M)$.
Recursively define 
${\mathcal A}_{k+1}^m(M)\subset {\mathcal A}_k^m(M)$
such that $A\in{\mathcal A}_k^m(M)$ belongs to ${\mathcal A}_{k+1}^m(M)$
if and only if $[A,D]=AD-DA\in{\mathcal A}_k^m(M)$ for every
smooth vector field $D$ on $M$.
Now we will use a variant of the commutator characterization
of pseudo-differential operators
(see e.g. \cite{Beals,
CoifmanMeyer,Co75,Dunau,Taylor97}), but we will
need the following Sobolev space version proved 
in \cite{Turunen00}, assuring that
the behaviour of commutators in Sobolev spaces
characterizes pseudo-differential operators:

\begin{thm}\label{THM:commutatorcharacterization}
A continuous linear operator $A:C^\infty(M)\to C^\infty(M)$
belongs to $\Psi^m(M)$ if and only if
$A\in\bigcap_{k=0}^\infty{\mathcal A}_k^m(M)$.
\end{thm}

In such characterization on a compact Lie group $M=G$,
it suffices to consider vector fields of the form
$D=M_\phi \partial_x$,
where $M_\phi f := \phi f$ is multiplication
by $\phi\in C^\infty(G)$,
and $\partial_x$ is left-invariant.
Notice that
$$
  [A,M_\phi \partial_x] = M_\phi\ [A,\partial_x] + [A,M_\phi]\ \partial_x,
$$
where $[A,M_\phi]f = A(\phi f)-\phi Af$.
Hence we need to consider compositions
$M_\phi A$, $A M_\phi$, $A\circ \partial_x$ 
and $\partial_x\circ A$.
First, we observe that 
\begin{eqnarray}
  \sigma_{M_\phi A}(x,\xi)
  & = & \phi(x)\ \sigma_A(x,\xi), \\
  \sigma_{A\circ\partial_x}(x,\xi)
  & = & \sigma_A(x,\xi)\ \sigma_{\partial_x}(x,\xi),\\
  \sigma_{\partial_x\circ A}(x,\xi)
  & = & \sigma_{\partial_x}(x,\xi)\ \sigma_A(x,\xi)
  +(\partial_x\sigma_A)(x,\xi),
  \label{partialAcomposition}
\end{eqnarray}
where $\sigma_{\partial_x}(x,\xi)$ is independent of $x\in G$.
Here \eqref{partialAcomposition} follows by the Leibnitz formula:
\begin{eqnarray*}
  \partial_x\circ A f(x)
  & = & \partial_x \sum_{[\xi]\in\widehat{G}} 
  {\rm dim}(\xi)\ {\rm Tr}\left(
      \xi(x)\ \sigma_A(x,\xi)\ \widehat{f}(\xi) \right) \\
    & = & \sum_{[\xi]\in\widehat{G}} {\rm dim}(\xi)\ {\rm Tr}\left(
      (\partial_x\xi)(x)\ \sigma_A(x,\xi)\ \widehat{f}(\xi) 
      \right) \\
    && + \sum_{[\xi]\in\widehat{G}} {\rm dim}(\xi)\ {\rm Tr}\left(
      \xi(x)\ \partial_x\sigma_A(x,\xi)\ \widehat{f}(\xi) \right).
\end{eqnarray*}
Next we claim that we have the fomula
\begin{equation}\label{compositionfunction}
  \sigma_{A M_\phi}(x,\xi) \sim \sum_{\alpha\geq 0} \frac{1}{\alpha!}
  \ \triangle_\xi^\alpha \sigma_A(x,\xi)\ \partial_x^{(\alpha)}\phi(x),
\end{equation}
where $\partial_x^{(\alpha)}$ are certain partial differential operators
of order $|\alpha|$.
This will follow from the following general composition formula:

\begin{thm}\label{THM:composition}
Let $m_1,m_2\in\Bbb R$ and $\rho>\delta\geq 0$.
Let $A,B:C^\infty(G)\to C^\infty(G)$ be continuous and linear,
their symbols satisfying
\begin{eqnarray*}
  \left\| \triangle_\xi^\alpha  \sigma_A(x,\xi) \right\|_{op}
  & \leq & C_{\alpha}\ 
  \langle\xi\rangle^{m_1-\rho|\alpha|},
  \\
  \left\| \partial_x^\beta \sigma_B(x,\xi) \right\|_{op}
  & \leq & C_{\beta}\ \langle\xi\rangle^{m_2+\delta|\beta|},
\end{eqnarray*}
for all multi-indices $\alpha$ and $\beta$,
uniformly in $x\in G$ and $[\xi]\in\widehat{G}$.
Then
\begin{equation}\label{compositiongeneral}
  \sigma_{AB}(x,\xi) \sim \sum_{\alpha\geq 0} \frac{1}{\alpha!}
  \ (\triangle_\xi^\alpha \sigma_A)(x,\xi)\ 
  \partial_x^{(\alpha)}\sigma_B(x,\xi),
\end{equation}
where the asymptotic expansion means that for every $N\in\N$
we have 
$$\n{\sigma_{AB}(x,\xi) -
   \sum_{|\alpha| < N} \frac{1}{\alpha!}
  \ (\triangle_\xi^\alpha \sigma_A)(x,\xi)\ 
  \partial_x^{(\alpha)}\sigma_B(x,\xi)}_{op}\leq
  C_N \jp{\xi}^{m_1+m_2-(\rho-\delta)N}.
$$
\end{thm}

\begin{proof}
First,
\begin{eqnarray*}
  ABf(x) & = & \int_G (Bf)(xz)\ R_A(x,z^{-1})\ {\rm d}z \\
  & = & \int_G \int_G f(xy^{-1})\ R_B(xz,yz)\ {\rm d}y
  \ R_A(x,z^{-1})\ {\rm d}z,
\end{eqnarray*}
where we use the standard distributional interpretation
of integrals.
Hence
\begin{eqnarray*}
  \sigma_{AB}(x,\xi)&=&\int_G R_{AB}(x,y)\ \xi(y)^\ast\ 
  {\rm d}y \\
  &=& \int_G \int_G R_A(x,z^{-1})\ \xi(z^{-1})^\ast\ 
  R_B(xz,yz)\ \xi(yz)^\ast
  \ {\rm d}z\ {\rm d}y \\ 
  & = & \sum_{|\alpha| < N} \frac{1}{\alpha!} \int_G \int_G
  R_A(x,z^{-1})\ q_\alpha(z^{-1})\ \xi(z^{-1})^\ast
  \ \partial_x^{(\alpha)} R_B(x,yz)\ \xi(yz)^\ast\ 
  {\rm d}z\ {\rm d}y \\
  && + \sum_{|\alpha| = N} \int_G \int_G
  R_A(x,z^{-1})\ q_\alpha(z^{-1})\ \xi(z^{-1})^\ast
  \ u_\alpha(x,yz)\ \xi^*(yz)\ {\rm d}z\ {\rm d}y \\
 & = & 
  \sum_{|\alpha| < N} \frac{1}{\alpha!}
  \ (\triangle_\xi^\alpha \sigma_A)(x,\xi)\ 
  \partial_x^{(\alpha)}\sigma_B(x,\xi) +
   \sum_{|\alpha|=N} 
  \ (\triangle_\xi^\alpha \sigma_A)(x,\xi)\ 
  \widehat{u_\alpha}(x,\xi).  
\end{eqnarray*}
Now the statement follows because we have
$\n{\widehat{u_\alpha}(x,\xi)}_{op}\leq C
\jp{\xi}^{m_1+\delta N}$ since $u_\alpha(x,y)$ is the
remainder in the Taylor expansion of $R_B(x,y)$ in $x$
only and so it satisfies similar estimates to those
of $\sigma_B$ with respect to $\xi$. 
This completes the proof.
\end{proof}
Before discussing symbol classes, let us complement
Theorem \ref{THM:composition} with a result about adjoint
operators:
\begin{thm}\label{THM:adjoint}
Let $m\in\Bbb R$ and $\rho>\delta\geq 0$.
Let $A:C^\infty(G)\to C^\infty(G)$ be continuous and linear,
with symbol $\sigma_A$ satisfying
\begin{equation}\label{EQ:adj-symbols-thm}
  \left\| \triangle_\xi^\alpha \partial_x^\beta 
  \sigma_A(x,\xi) \right\|_{op}
   \leq  C_{\alpha}\ 
  \langle\xi\rangle^{m-\rho|\alpha|+\delta|\beta|},
\end{equation}
for all multi-indices $\alpha$,
uniformly in $x\in G$ and $[\xi]\in\widehat{G}$.
Then the symbol of $A^*$ is
\begin{equation}\label{compositiongeneral2}
  \sigma_{A^*}(x,\xi) \sim \sum_{\alpha\geq 0} \frac{1}{\alpha!}
  \ \triangle_\xi^\alpha\partial_x^{(\alpha)} 
  \sigma_A(x,\xi)^*,
\end{equation}
where the asymptotic expansion means that for every $N\in\N$
we have 
$$\n{\triangle_\xi^\gamma\partial_x^\beta\p{\sigma_{A}(x,\xi) -
   \sum_{|\alpha| < N} \frac{1}{\alpha!}
 \ \triangle_\xi^\alpha\partial_x^{(\alpha)} 
  \sigma_A(x,\xi)^*}}_{op}\leq
  C_N \jp{\xi}^{m-(\rho-\delta)N-\rho|\gamma|+\delta|\beta|}.
$$
\end{thm}
\begin{rem}
We note that if we impose conditions of the type
\eqref{EQ:adj-symbols-thm} on both symbols
$\sigma_A, \sigma_B$ in Theorem \ref{THM:composition},
we also get the asymptotic expansion 
\eqref{compositiongeneral} with the remainder estimate as in 
Theorem \ref{THM:adjoint}.
\end{rem}
\begin{proof}[Proof of Theorem \ref{THM:adjoint}.]
First we observe that writing
$A^*g(y)=\int_G g(x) R_{A^*}(y,x^{-1}y)\drm x$, we
get the relation
$R_{A^*}(y,x^{-1}y)=\overline{R_{A}(x,y^{-1}x)}$ between kernels,
which means that
$R_{A^*}(x,v)=\overline{R_{A}(x v^{-1},v^{-1})}$.
From this we find
\begin{eqnarray*}
  \sigma_{A^*}(x,\xi)&=&\int_G R_{A^*}(x,v)\ \xi(v)^\ast\ 
  {\rm d}v \\
  &=& \int_G \overline{R_{A}(x v^{-1},v^{-1})}\ 
  \xi(v)^\ast\ 
  {\rm d}v \\ 
  & = & \sum_{|\alpha| < N} \frac{1}{\alpha!} \int_G 
   q_\alpha(v)\ \partial_x^{(\alpha)}
   \overline{R_A(x,v^{-1})}\  \xi(v)^\ast
   {\rm d}v  + \Rcal_N(x,\xi)  \\
 & = & 
 \sum_{|\alpha| < N} \frac{1}{\alpha!}
 \ \triangle_\xi^\alpha\partial_x^{(\alpha)} 
  \sigma_A(x,\xi)^* + \Rcal_N(x,\xi), 
\end{eqnarray*}
where the last formula for the asymptotic expansion
follows in view of
$$
\sigma_A(x,\xi)^*=\p{\int_G R_A(x,v)\ \xi^*(v) \drm v}^*
 = \int_G \overline{R_A(x,v^{-1})}\ \xi^*(v) \drm v,
$$
and estimate for the remainder $\Rcal_N(x,\xi)$ follows
by an argument similar to that in the proof of
Theorem \ref{THM:composition}.
\end{proof} 

On the way to characterize
the usual H\"ormander's classes $\Psi^m(G)$ in
Theorem \ref{THM:Gcommutatorcharacterization},
we need some properties concerning
symbols of pseudo-differential operators.

\begin{lem}\label{LEM:freezing}
Let $A\in\Psi^m(G)$. 
Then there exists a constant $C<\infty$
such that
$$
  \left\|\sigma_A(x,\xi)\right\|_{op}\leq C\langle\xi\rangle^m
$$
for all $x\in G$ and $\xi\in\Rep(G)$.
Also, if $u\in G$ and if $B$ is an 
operator with symbol
 $\sigma_{B}(x,\xi)=\sigma_A(u,\xi)$, then
$B\in\Psi^m(G)$.
\end{lem}
\begin{proof}
First, $B\in\Psi^m(G)$
follows from the local theory of 
pseudo-differential operators, 
by studying 
$B f(x) = \int_G K_A(u,u x^{-1} y)\ f(y)\ {\rm d}y$.
Hence the right-convolution operator $B$
is bounded from $H^s(M)$ to $H^{s-m}(G)$,
implying $\|\sigma_A(u,\xi)\|\leq C \langle\xi\rangle^m$.
\end{proof}

\begin{lem}\label{LEM:pointwisesymbols}
Let $A\in\Psi^m(G)$. Then
$Op(\triangle_\xi^\alpha\partial_x^\beta\sigma_A)\in\Psi^{m-|\alpha|}(G)$
for all $\alpha,\beta$.
\end{lem}

\begin{proof}
First, given $A\in\Psi^m(G)$, let us define
$\sigma_B(x,\xi)=\triangle_\xi^\alpha\partial_x^\beta\sigma_A(x,\xi)$.
We must show that $B\in\Psi^{m-|\alpha|}(G)$.
If here $|\beta|=0$, we obtain
$$
  Bf(x)  =  \int_G f(xy^{-1})\ q_\alpha(y)\ R_A(x,y)\ 
  {\rm d}y 
   =  \int_G q_\alpha(y^{-1}x)\ K_A(x,y)\ f(y)\ {\rm d}y.
$$
Moving to local coordinates, we need to study
$$
  \tilde{B}f(x) = \int_{\Bbb R^n}
  \phi(x,y)\ K_{\tilde{A}}(x,y)\ f(y)\ {\rm d}y,
$$
where $\tilde{A}\in\Psi^m(\Bbb R^n)$
with $\phi\in C^\infty(\Bbb R^n\times\Bbb R^n)$,
the kernel $K_{\tilde{A}}$ being compactly supported.
Let us calculate the symbol of $\tilde{B}$:
\begin{eqnarray*}
  \sigma_{\tilde{B}}(x,\xi)
  & = & \int_{\Bbb R^n}\ {\rm e}^{{\rm i}2\pi(y-x)\cdot\xi}
  \ \phi(x,y)\ K_{\tilde{A}}(x,y)\ {\rm d}y \\
  & \sim & \sum_{\gamma\geq 0} \frac{1}{\gamma!}
  \ \left.\partial_z^\gamma \phi(x,z)\right|_{z=x}
  \int_G {\rm e}^{{\rm i}2\pi(y-x)\cdot\xi}\ (y-x)^{\gamma}
  \ K_{\tilde{A}}(x,y)\ {\rm d}y \\
  & = & \sum_{\gamma\geq 0} \frac{1}{\gamma!}
  \ \left.\partial_z^\gamma \phi(x,z)\right|_{y=x}
  D_\xi^\gamma\sigma_{\tilde{A}}(x,\xi).
\end{eqnarray*}
This shows that $\tilde{B}\in\Psi^m(\Bbb R^n)$.
We obtain $Op(\triangle_\xi^\alpha\sigma_A)\in\Psi^{m-|\alpha|}(G)$
if $A\in\Psi^m(G)$.

Next we show that $B=Op(\partial_x^\beta\sigma_A)\in\Psi^m(G)$.
We may assume that $|\beta|=1$.
Left-invariant vector field $\partial_x^\beta$ is a linear combination
of terms of the type $c(x) D_x$,
where $c\in C^\infty(G)$ and $D_x$ is right-invariant.
By the previous considerations on $\tilde{B}$,
we may remove $c(x)$ here, and consider only
$C=Op(D_x\sigma_A)$.
Since $R_A(x,y) = K_A(x,xy^{-1})$, we get
$$
  D_x R_A(x,y) = \left.(D_x + D_z) K_A(x,z)\right|_{z=xy^{-1}},
$$
leading to
$$
  Cf(x)  =  \int_G f(xy^{-1})\ D_x R_A(x,y)\ {\rm d}y 
   =  \int_G f(y)\ (D_x+D_y) K_A(x,y)\ {\rm d}y.
$$
Thus, we study local operators of the form
\begin{eqnarray*}
  \tilde{C}f(x) & = & \int_{\Bbb R^n} f(y)
  \left( \phi(x,y)\partial_x^{\beta} + \psi(x,y)
  \partial_y^\beta \right)
  K_{\tilde{A}}(x,y)\ {\rm d}y,
\end{eqnarray*}
where the kernel of $\tilde{A}\in\Psi^m(\Bbb R^n)$ has compact support,
$\phi,\psi\in C^\infty(\Bbb R^n\times\Bbb R^n)$,
and $\phi(x,x)=\psi(x,x)$ for every $x\in\Bbb R^n$.
Let $\tilde{C}=\tilde{D}+\tilde{E}$, where
\begin{eqnarray*}
  \tilde{D}f(x) & = & \int_{\Bbb R^n} f(y)
  \ \phi(x,y) \left(\partial_x^{\beta} + \partial_y^\beta \right)
  K_{\tilde{A}}(x,y)\ {\rm d}y, \\
  \tilde{E}f(x)
  & = & \int_{\Bbb R^n} f(y)\ \left(\psi(x,y)-\phi(x,y)\right)
  \ \partial_y^\beta K_{\tilde{A}}(x,y)\ {\rm d}y.
\end{eqnarray*}
By the above considerations about $\tilde{B}$,
we may assume that $\phi(x,y)\equiv 1$ here, and obtain
$\sigma_{\tilde{D}}(x,\xi) = \partial_x^\beta\sigma_{\tilde{A}}(x,\xi)$.
Thus $\tilde{D}\in\Psi^m(\Bbb R^n)$.
Moreover,
\begin{eqnarray*}
  \tilde{E}f(x)
  & \sim & \sum_{\gamma\geq 0}\frac{1}{\gamma!}
  \ \left. \partial_z^\gamma \left(\psi(x,z)-\phi(x,z)\right)\right|_{z=x}
  \int_{\Bbb R^n} f(y)
  \ (y-x)^\gamma\ \partial_y^\beta K_{\tilde{A}}(x,y) \ {\rm d}y,
\end{eqnarray*}
yielding
\begin{eqnarray*}
  \sigma_{\tilde{E}}(x,\xi)
  & \sim & \sum_{\gamma\geq 0} c_\gamma(x)\ \partial_\xi^\gamma
  \left(\xi^\beta\ \sigma_{\tilde{A}}(x,\xi)\right)
\end{eqnarray*}
for some functions $c_\gamma\in C^\infty(\Bbb R^n)$
for which $c_0(x)\equiv 0$.
Since $|\beta|=1$, this shows that $\tilde{E}\in\Psi^m(\Bbb R^n)$.
Thus $Op(\partial_x^\beta\sigma_A)\in\Psi^m(G)$
if $A\in\Psi^m(G)$.
\end{proof}

\begin{lem}\label{LEM:symbolproductderivative}
Let $A\in\Psi^m(G)$ and let $D:C^\infty(G)\to C^\infty(G)$
be a smooth vector field.
Then $Op(\sigma_A \sigma_D)\in\Psi^{m+1}(G)$
and $Op([\sigma_A,\sigma_D])\in\Psi^m(G)$.
\end{lem}

\begin{proof}
For simplicity, we may assume that $D=M_\phi\partial_x$,
where $\partial_x$ is left-invariant
and $\phi\in C^\infty(G)$.
Now
$$
  \sigma_A(x,\xi)\ \sigma_D(x,\xi)
   =  \phi(x)\ \sigma_A(x,\xi)\ \sigma_{\partial_x}(\xi) 
   =  \sigma_{M_\phi A\circ\partial_x}(x,\xi),
$$
and it is well-known that $M_\phi A\circ\partial_x\in\Psi^{m+1}(G)$.
Thus $Op(\sigma_A \sigma_D)\in\Psi^{m+1}(G)$.
Next,
\begin{eqnarray*}
  \sigma_D(x,\xi)\ \sigma_A(x,\xi)
  & = & \phi(x)\ \sigma_{\partial_x}(\xi)\ \sigma_A(x,\xi) \\
  & \stackrel{\eqref{partialAcomposition}}{=} &
  \phi(x)\left( \sigma_{\partial_x\circ A}(x,\xi)
    - (\partial_x\sigma_A)(x,\xi) \right) \\
  & = & \sigma_{M_\phi\circ\partial_x\circ A}(x,\xi)
  - \phi(x)\ (\partial_x\sigma_A)(x,\xi).
\end{eqnarray*}
From this we see that
\begin{eqnarray*}
  Op([\sigma_A,\sigma_D]) & = & M_\phi A\circ\partial_x
  - M_\phi \partial_x\circ A + M_\phi\ Op(\partial_x\sigma_A) \\
  & = & M_\phi [A,\partial_x]+ M_\phi\ Op(\partial_x\sigma_A).
\end{eqnarray*}
Here $Op(\partial_x\sigma_A)\in\Psi^m(G)$
by Lemma~\ref{LEM:pointwisesymbols}.
Hence $Op([\sigma_A,\sigma_D])$
belongs to $\Psi^m(G)$ by the known properties
of pseudo-differential operators.
\end{proof}

Finally, let us prove the Sobolev space boundedness
of pseudo-differential operators given in
Theorem~\ref{THM:su2-Sobolev}.

\begin{proof}[Proof of Theorem~\ref{THM:su2-Sobolev}.]
Observing the continuous
mapping $\Xi^s:H^s(G)\to L^2(G)$, we have to 
prove that operator $\Xi^{s-\mu}\circ A\circ \Xi^{-s}$
is bounded from $L^2(G)$ to $L^2(G)$. Let us denote
$B=A\circ \Xi^{-s}$, so that the symbol of $B$ satisfies
$\sigma_B(x,\xi)=\jp{\xi}^{-s}\sigma_A(x,\xi)$ for
all $x\in G$ and $\xi\in\Rep(G)$, where $\langle\xi\rangle$ is defined
in \eqref{EQ:su2-jp-bracket}. 
Since $\Xi^{s-\mu}\in\Psi^{s-\mu}(G)$, 
by \eqref{EQ:su2-l2-norm-mult} and 
Lemma \ref{LEM:pointwisesymbols}
its symbol satisfies 
\begin{equation}\label{EQ:su2-Sob-op}
  \|\triangle_\xi^\alpha\sigma_{\Xi^{s-\mu}}(x,\xi)\|_{op}
  \leq C'_\alpha\ \jp{\xi}^{s-\mu-|\alpha|}.
\end{equation} 
Now we can observe that
the asymptotic formula in Theorem~\ref{THM:composition}
works for the composition $\Xi^{s-\mu}\circ B$ in view
of \eqref{EQ:su2-Sob-op},
and we obtain
\begin{eqnarray*}
  \partial_x^{\beta} \sigma_{\Xi^{s-\mu}\circ B}(x,\xi)
  & \sim & \sum_{\alpha\geq 0} \frac{1}{\alpha!}
  \left( \triangle_\xi^\alpha \sigma_{\Xi^{s-\mu}}(x,\xi) \right)
  \ \langle\xi\rangle^{-s}
  \ \partial_x^{(\alpha)} \partial_x^\beta \sigma_A(x,\xi).
\end{eqnarray*}
It follows that
$$
  \left\| \partial_x^\beta
  \sigma_{\Xi^{s-\mu}\circ B}(x,\xi) \right\|_{op}
  \leq C''_\beta,
$$
so that $\Xi^{s-\mu}\circ B$ is
bounded on $L^2(G)$ by Theorem~\ref{THM:su2-L2}.
This completes the proof.
\end{proof}

\section{Symbol classes on compact Lie groups}
\label{SEC:symbolclasses}

The goal of this section is to describe
the pseudo-differential symbol inequalities on compact Lie groups
that yield H\"ormander's classes $\Psi^m(G)$.
Combined with
asymptotic expansion~\eqref{compositiongeneral}
for composing operators,
Theorem~\ref{THM:commutatorcharacterization}
motivates defining the following symbol classes
$\Sigma^m(G)=\bigcap_{k=0}^\infty \Sigma_k^m(G)$, that 
we will show to characterize H\"ormander's class
$\Psi^m(G)$.

\begin{defn}\label{symbolinequalities}
Let $m\in\R$.
We denote $\sigma_A\in \Sigma_0^m(G)$ if 
\begin{equation} \label{cond5-G}
  {\rm sing\ supp}\left(y\mapsto R_A(x,y)\right)
   \subset \{e\}
\end{equation} 
and if
\begin{equation}\label{norminequality}
  \| \triangle_\xi^\alpha \partial_x^\beta \sigma_A(x,\xi) \|_{op}
        \leq C_{A\alpha\beta m}\  
        \langle \xi \rangle^{m-|\alpha|},
\end{equation}
for all $x\in G$, all multi-indices
$\alpha,\beta$, and all $\xi\in\Rep(G)$, 
where $\jp{\xi}$ is defined in \eqref{EQ:su2-jp-bracket}.
Then we say that
$\sigma_A\in \Sigma_{k+1}^m(G)$ if and only if
\begin{eqnarray}
  \sigma_A
  & \in &
        \Sigma_k^m(G),\\
  \sigma_{\partial_j}\sigma_A - \sigma_A\sigma_{\partial_j}
  & \in &
        \Sigma_k^m(G),\\
\label{su2-symbols-G-cond4}
  (\triangle_\xi^\gamma\sigma_A)\ \sigma_{\partial_j}
  & \in &
        \Sigma_k^{m+1-|\gamma|}(G),
\end{eqnarray}
for all $|\gamma|>0$ and $1\leq j\leq{\rm dim}(G)$.
Let
\begin{equation*}
  \Sigma^m(G) = \bigcap_{k=0}^\infty \Sigma_k^m(G).
\end{equation*}
Let us denote $A\in Op\Sigma^m(G)$ if and only if $\sigma_A\in\Sigma^m(G)$.
\end{defn}

\begin{thm}\label{THM:Gcommutatorcharacterization}
Let $G$ be a compact Lie group and let $m\in\R$.
Then $A\in \Psi^m(G)$ if and only if $\sigma_A\in \Sigma^m(G)$,
i.e. $Op\Sigma^m(G)=\Psi^m(G)$.
\end{thm}

\begin{proof}
First, applying Theorem~\ref{THM:composition} to 
$\sigma_A\in \Sigma_{k+1}^m(G)$,
we notice that $[A,D]\in Op\Sigma_k^m(G)$
for any smooth vector field $D:C^\infty(G)\to C^\infty(G)$.
Consequently, if here $A\in Op\Sigma^m(G)$
then also $[A,D]\in Op\Sigma^m(G)$.
By Remark~\ref{REM:SobolevL2boundedness},
$Op\Sigma^m(G)\subset{\mathcal L}(H^m(G),L^2(G))$.
Hence Theorem~\ref{THM:commutatorcharacterization} implies
$Op\Sigma^m(G)\subset\Psi^m(G)$.

Conversely, we have to show that $\Psi^m(G)\subset Op\Sigma^m(G)$.
This follows by Lemma~\ref{LEM:freezing},
and Lemmas~\ref{LEM:pointwisesymbols} and \ref{LEM:symbolproductderivative}.
More precisely, let $A\in\Psi^m(G)$.
Then we have 
\begin{eqnarray*}
  Op(\triangle_\xi^\alpha\partial_x^\beta\sigma_A)
  & \in & \Psi^{m-|\alpha|}(G), \\
  Op\left([\sigma_{\partial_j},\sigma_A]\right)
  & \in & \Psi^m(G), \\
  Op\left((\triangle_\xi^\gamma\sigma_A)\sigma_{\partial_j}\right)
  & \in & \Psi^{m+1-|\gamma|}(G).
\end{eqnarray*}
Moreover, $\left\|\sigma_A(x,\xi)\right\|\leq C\langle\xi\rangle^m$
by Lemma~\ref{LEM:freezing},
and the singular support $y\mapsto R_A(x,y)$
is contained in $\{e\}\subset G$.
This completes the proof.
\end{proof}

\begin{cor}\label{COR:psdosymbolinvariances}
The set $\Sigma^m(G)$ is invariant under $x$-freezings,
$x$-translations and $\xi$-conjugations.
More precisely, if $(x,\xi)\mapsto\sigma_A(x,\xi)$ belongs to 
$\Sigma^m(G)$
and $u\in G$ then also the following symbols belong to 
$\Sigma^m(G)$:
\begin{eqnarray}
  (x,\xi) & \mapsto & \sigma_A(u,\xi), \label{freezing} \\
  (x,\xi) & \mapsto & \sigma_A(ux,\xi), \label{Ltrans} \\
  (x,\xi) & \mapsto & \sigma_A(xu,\xi), \label{Rtrans} \\
  (x,\xi) & \mapsto & \xi(u)^\ast\ \sigma_A(x,\xi)\ \xi(u). 
  \label{Tconj}
\end{eqnarray}
\end{cor}

\begin{proof}
The symbol classes $\Sigma^m(G)$ are defined by conditions
\eqref{cond5-G}-\eqref{su2-symbols-G-cond4},
which are checked for points $x\in G$ fixed
(with constants uniform in $x$).
Therefore it follows that $\Sigma^m(G)$ is invariant
under the $x$-freezing (\ref{freezing}),
and under the left and right $x$-translations 
\eqref{Ltrans},\eqref{Rtrans}.
The $x$-freezing property (\ref{freezing}) 
would have followed also from
Lemma~\ref{LEM:freezing} and Theorem~\ref{THM:Gcommutatorcharacterization}.
From the general theory of pseudo-differential operators
it follows that $A\in\Psi^m(G)$
if and only if the $\phi$-pullback $A_\phi$ belongs to the same class
$\Psi^m(G)$,
where
$
  A_\phi f = A(f\circ\phi)\circ\phi^{-1}.
$
This, combined with the $x$-translation invariances
and Proposition~\ref{PROP:conjugatedrightsymbols},
implies the conjugation invariance in (\ref{Tconj}).
\end{proof}
From Theorem \ref{THM:Gcommutatorcharacterization} and
Lemma \ref{LEM:pointwisesymbols} we also obtain:
\begin{cor}\label{COR:G-sigma-diff}
If $\sigma_A\in \Sigma^m(G)$ then
$\triangle_\xi^\alpha \partial_x^\beta \sigma_A
\in \Sigma^{m-|\alpha|}(G)$.
\end{cor}

\section{Symbol classes on $\SU2$}

Let us now turn to the analysis on $\SU2$. In this section
we derive a much simpler symbolic characterization of 
pseudo-differential operators on $\SU2$ than the one 
given in Definition~\ref{symbolinequalities}.
First we summarize the approach in the case of
$\SU2$ also simplifying the notation in this case.

By the Peter-Weyl theorem
$\{\sqrt{2l+1}\ {t^l_{nm}}:\ l\in\frac12\N_0,
        \ -l\leq m,n\leq l,\ l-m,l-n\in\Bbb Z\}$
is an orthonormal basis for $L^2(\SU2)$, where 
$t^l$ were defined in Section~\ref{SEC:preliminaries},
and thus $f\in C^\infty(\SU2)$ has a Fourier series representation
\begin{equation*}
  f = \sum_{l\in\frac12 \N_0}(2l+1) \sum_m \sum_n
        \widehat{f}(l)_{mn}\ {t^l_{nm}},
\end{equation*}
where the Fourier coefficients are computed by
\begin{equation*}
  \widehat{f}(l)_{mn} := 
       \int_\SU2 f(g)\ \overline{t^l_{nm}(g)}\ {\rm d}g =
       \langle f,{t^l_{nm}}\rangle_{L^2(\SU2)},
\end{equation*}
so that $\widehat{f}(l)\in\C^{(2l+1)\times(2l+1)}$.
We recall that in the case of $\SU2$, we simplify the notation
writing $\widehat{f}(l)$ instead of $\widehat{f}(t^l)$, etc.

Let $A:C^\infty(\SU2)\to C^\infty(\SU2)$ 
be a continuous linear operator
and let $R_A\in{\mathcal D}'(\SU2\times \SU2)$ be
its right-convolution kernel, i.e.
$$
  Af(x) = \int_\SU2 f(y)\ R_A(x,y^{-1}x)\ {\rm d}y
  = (f\ast R_A(x,\cdot))(x)
$$
in the sense of distributions.
According to Definition~\ref{DEF:su2-symbols-on-G},
by the symbol of $A$
we mean the sequence of matrix-valued mappings
$$
  (x\mapsto\sigma_A(x,l)):\SU2\to\Bbb C^{(2l+1)\times(2l+1)},
$$
where $2l\in\Bbb N_0$, obtained from
$$
  \sigma_A(x,l)_{mn} = \int_\SU2 R_A(x,y)\ \overline{t^l_{nm}(y)}
  \ {\rm d}y.
$$
That is, $\sigma_A(x,l)$ is the $l^{\rm th}$ Fourier coefficient
of the function $y\mapsto R_A(x,y)$.
Then by Theorem~\ref{THM:su2-symbol} we have
\begin{eqnarray*}
  Af(x)
  & = & \sum_{l} (2l+1)\ {\rm Tr}\left( 
  t^l(x) \ \sigma_A(x,l) \ \widehat{f}(l)
     \right) \\
  & = & \sum_{l} (2l+1)\ \sum_{m,n} 
   t^l(x)_{nm}\
  \left( \sum_k  \sigma_A(x,l)_{mk}
    \ \widehat{f}(l)_{kn} \right).
\end{eqnarray*}
Alternatively, by Theorem~\ref{THM:su2-symbol2} we have
\begin{equation*}
  \sigma_A(x,l) = t^l(x)^\ast\
   \left( A t^l\right)(x),
\end{equation*}
that is
\begin{equation*}
  \sigma_A(x,l)_{mn} 
  =  \sum_k \overline{t^l_{km}(x)}(At^l_{kn})(x).
\end{equation*}

In the case of $\SU2$, quantity 
$\jp{t^l}$ for $\xi=t^l$ in \eqref{EQ:su2-jp-bracket}
can be calculated as
$$
 \jp{t^l}=(1+\lambda_{[t^l]})^{1/2}=
 (1+l(l+1))^{1/2} 
$$
in view of Theorem~\ref{THM:tlmnpartialderivatives},
and Definition \ref{symbolinequalities} becomes:


\begin{defn}
We write that symbol $\sigma_A\in \Sigma_0^m(\SU2)$ if 
\begin{equation} \label{cond5}
  {\rm sing\ supp}\left(y\mapsto R_A(x,y)\right)
   \subset \{e\}
\end{equation} 
and if
\begin{equation}\label{cond0}
  \left\| \triangle_l^\alpha \partial_x^\beta \sigma_A(x,l)
  \right\|_{\Bbb C^{2l+1}\to\Bbb C^{2l+1}}
        \leq C_{A\alpha\beta m}\  (1+l)^{m-|\alpha|}
\end{equation}
for all $x\in G$, all multi-indices $\alpha,\beta$,
 and $l\in\frac{1}{2}\Bbb N_0$.
Here $\triangle_l^\alpha=\triangle_0^{\alpha_1}\triangle_+^{\alpha_2}
\triangle_-^{\alpha_3}$ and
$\partial_x^\beta=\partial_0^{\beta_1}\partial_+^{\beta_2}
\partial_-^{\beta_3}$.
Moreover, $\sigma_A\in \Sigma_{k+1}^m(\SU2)$ if and only if
\begin{eqnarray}
  \sigma_A
  & \in &
        \Sigma_k^m(\SU2),\label{cond1} \\
  {}[\sigma_{\partial_j},\sigma_A]
  = \sigma_{\partial_j}\sigma_A - \sigma_A\sigma_{\partial_j}
  & \in &
        \Sigma_k^m(\SU2),\label{cond2} \\
  (\triangle_l^\gamma\sigma_A)\ \sigma_{\partial_j}
  & \in &
        \Sigma_k^{m+1-|\gamma|}(\SU2),\label{cond4}
\end{eqnarray}
for all $|\gamma|>0$ and $j\in\{0,+,-\}$.
Let
\begin{equation*}
  \Sigma^m(\SU2) = \bigcap_{k=0}^\infty \Sigma_k^m(\SU2),
\end{equation*}
so that by Theorem~\ref{THM:Gcommutatorcharacterization}
we have $Op \Sigma^m(\SU2)=\Psi^m(\SU2)$.
\end{defn}

\begin{rem}
We would like to provide a more direct definition for $\Sigma^m(\SU2)$,
without resorting to classes $\Sigma_k^m(\SU2)$.
Condition (\ref{cond0})
is just an analogy of the usual symbol inequalities.
Conditions \eqref{cond5} and \eqref{cond1} are straightforward.
We may have difficulties with differences $\triangle_l^\alpha$,
but derivatives $\partial_x^\beta$ do not cause problems;
if we want, we may assume that the symbols are constant in $x$.
By the definition of operators $\triangle_l^\alpha$ and
$\partial_x^\beta$ we also have the following properties:
$$
  \triangle_l^\alpha\partial_x^\beta \sigma_A(x,l)
  = \partial_x^\beta\triangle_l^\alpha \sigma_A(x,l),
$$
$$
  \partial_j\left(\sigma_A(x,l)\ \sigma_B(x,l) \right)
  = \left( \partial_j \sigma_A(x,l) \right) \sigma_B(x,l)
  + \sigma_A(x,l)\ \partial_j \sigma_A(x,l),
$$
$$
  \partial_y^\beta \left(
    \sigma_A(x,l)\ \sigma_B(y,l)\ \sigma_C(z,l) \right)
  = \sigma_A(x,l) \left(\partial_y^\beta\sigma_B(y,l)\right)
  \sigma_C(z,l).
$$
\end{rem}

We now give another, simpler characterization
of pseudo-differential operators.

\begin{defn}[Symbol classes on $\SU2$]
For $u\in{\rm SU}(2)$,
denote $A_u f:= A(f\circ\phi)\circ\phi^{-1}$,
where $\phi(x)=xu$;
then (by Proposition \ref{PROP:conjugatedrightsymbols})
\begin{eqnarray*}
  R_{A_u}(x,y) & = & R_A(xu^{-1},uyu^{-1}), \\
  \sigma_{A_u}(x,l) & = & t^l(u)^\ast\ \sigma_A(xu^{-1},l)\ t^l(u).
\end{eqnarray*}
The symbol class $S^m({\rm SU}(2))$
consists of the symbols $\sigma_A$
of those operators $A\in{\mathcal L}(C^\infty({\rm SU}(2)))$
for which $(y\mapsto R_A(x,y))\subset\{e\}$ and
for which 
\begin{equation}\label{Sigmainequalities}
  \left| \triangle_l^\alpha \partial_x^\beta 
  \sigma_{A_u}(x,l)_{ij}\right|
  \leq C_{A\alpha\beta m N}\ \langle i-j\rangle^{-N}
  \ (1+l)^{m-|\alpha|}
\end{equation}
uniformly in $x,u\in{\rm SU}(2)$, for every $N\geq 0$,
all $l\in\frac12\N_0$,
every multi-indices $\alpha,\beta\in\Bbb N_0^3$,
and for all matrix column/row numbers $i,j$.
Thus, the constant in \eqref{Sigmainequalities} may depend
on $A, \alpha, \beta, m$ and $N$, but not on
$x, u, l, i, j$.
\end{defn}

We now formulate the main theorem of this section:
\begin{thm}\label{THM:simplecharacterization}
Operator $A\in{\mathcal L}(C^\infty({\rm SU}(2)))$
belongs to $\Psi^m({\rm SU}(2))$
if and only if $\sigma_A\in S^m({\rm SU}(2))$.
Moreover, we have  the equality of symbol classes
 $S^m({\rm SU}(2))=\Sigma^m({\rm SU}(2))$.
\end{thm}
In fact, we need to prove only the equality of symbol classes
 $S^m({\rm SU}(2))=\Sigma^m({\rm SU}(2))$, from which
the first part of the theorem would follow
by Theorem \ref{THM:Gcommutatorcharacterization}.
In the process of proving this equality, we 
establish a number of auxiliary results.

\begin{rem}
By Corollary \ref{COR:G-sigma-diff},
if $\sigma_A\in \Sigma^m({\rm SU}(2))$ then
$\triangle_l^\gamma \partial_x^\delta \sigma_A
\in \Sigma^{m-|\delta|}({\rm SU}(2))$.
Let us show the analogous result for $S^m({\rm SU}(2))$.
\end{rem}

\begin{lem}\label{LEM:tildeSdiffbehavior}
If $\sigma_A\in S^m({\rm SU}(2))$ then
$\sigma_B = \triangle_l^\gamma \partial_x^\delta \sigma_A
\in S^{m-|\gamma|}({\rm SU}(2))$.
\end{lem}

\begin{proof}
First, let $|\gamma|=1$.
Then $\triangle_l^\gamma\widehat{f}(l)=\widehat{qf}(l)$
for some
$$
  q\in Pol_1({\rm SU}(2))
  :={\rm span}\left\{t^{1/2}_{ij}:\ i,j\in\{-1/2,+1/2\}\right\}
$$
for which $q(e)=0$.
Let $r(y):=q(u y u^{-1})$.
Then $r\in Pol_1({\rm SU}(2))$, because
$$
  t^{1/2}_{ij}(uyu^{-1}) = \sum_{k,m}
  t^{1/2}_{ik}(u)\ t^{1/2}_{km}(y)\ t^{1/2}_{mj}(u^{-1}).
$$
Moreover, we have $r(e)=0$.
Hence $\widehat{f}(l)\mapsto\widehat{rf}(l)$
is a linear combination of difference operators
$\triangle_0,\triangle_+,\triangle_-$
because $\left\{f\in Pol^1({\rm SU}(2)):\ f(e)=0\right\}$
is a three-dimensional vector space spanned by $q_0,q_+,q_-$.
Now let $\gamma\in\Bbb N_0^3$ and
$\sigma_B = \triangle_l^\gamma \partial_x^\delta \sigma_A$.
We have
\begin{eqnarray*}
  \triangle_l^\alpha \partial_x^\beta\sigma_{B_u}(x,l)
  & = & \triangle_l^\alpha \partial_x^\beta
  \left( t^l(u)^\ast\ \sigma_B(xu^{-1},l)\ t^l(u) \right) \\
  & = & \triangle_l^\alpha \partial_x^\beta
  \left( t^l(u)^\ast
  \left(\triangle_l^\gamma \partial_x^\delta 
  \sigma_A(xu^{-1},l)\right)
  t^l(u) \right) \\
  & = & \sum_{|\gamma'|=|\gamma|}
  \lambda_{u,\gamma'}
  \ \triangle_l^{\alpha+\gamma'} \partial_x^{\beta+\delta}
  \ \sigma_{A_u}(x,l),
\end{eqnarray*}
for some scalars $\lambda_{u,\gamma'}\in\Bbb C$ depending only on
$u\in{\rm SU}(2)$ and multi-indices $\gamma'\in\Bbb N_0^3$.
\end{proof}

\begin{rem}
Let $D$ be a left-invariant vector field on ${\rm SU}(2)$.
From the very definition of the symbol classes
$\Sigma^m({\rm SU}(2))=\bigcap_{k=0}^\infty \Sigma_k^m({\rm SU}(2))$,
it is evident that
$[\sigma_D,\sigma_A]\in \Sigma^m({\rm SU}(2))$
if $\sigma_A\in \Sigma^m({\rm SU}(2))$.
We shall next prove the similar invariance for $S^m({\rm SU}(2))$.
\end{rem}

\begin{lem}\label{LEM:tildeSproductbehavior}
Let $D$ be a left-invariant vector field on ${\rm SU}(2)$.
Let $\sigma_A\in S^m({\rm SU}(2))$.
Then $[\sigma_D,\sigma_A]\in S^m({\rm SU}(2))$ and
$\sigma_A\ \sigma_D\in S^{m+1}({\rm SU}(2))$.
\end{lem}

\begin{proof}
For $D\in {\mathfrak su(2)}$ we write $D=\irm E$,
so that $E\in \irm\ {\mathfrak su(2)}$.
By Proposition~\ref{PROP:SU2-vector-fields}
there is some $u\in{\rm SU}(2)$ such that
$\sigma_E(l) = t^l(u)^\ast\ \sigma_{\partial_0}(l)\ t^l(u)$.
Now, we have
\begin{eqnarray*}
  [\sigma_E,\sigma_A](l)
  & = & t^l(u)^\ast\ \left[\sigma_{\partial_0}(l),
    t^l(u)\ \sigma_A(x,l)\ t^l(u)^\ast \right]\ t^l(u) \\
  & = & \left[\sigma_{\partial_0},\sigma_{A_{u^{-1}}}\right]_u (l).
\end{eqnarray*}
Next, notice that $S^m({\rm SU}(2))$ is invariant under
the mappings $\sigma_B\mapsto\sigma_{B_u}$ and
$\sigma_B\mapsto[\sigma_{\partial_0},\sigma_B]$;
here $[\sigma_{\partial_0},\sigma_B](l)_{ij} = 
(i-j)\ \sigma_B(l)_{ij}$.
Finally,
\begin{eqnarray*}
  \sigma_A(x,l)\ \sigma_E(l)
  & = &
  t^l(u)^\ast\ t^l(u)\ \sigma_A(x,l)\ t^l(u)^\ast\ 
  \sigma_{\partial_0}(l)
  \ t^l(u) \\
  & = & \left( \sigma_{A_{u^{-1}}}(x,l)\ \sigma_{\partial_0}(l)
  \right)_u.
\end{eqnarray*}
Just like in the first part of the proof,
we see that $\sigma_A\ \sigma_D$
belongs to $S^{m+1}({\rm SU}(2))$ since
$\sigma_B\ \sigma_{\partial_0}\in S^{m+1}({\rm SU}(2))$
if $\sigma_B\in S^m({\rm SU}(2))$,
by Theorem~\ref{THM:su2-Leibnitz}.
\end{proof}

\begin{proof}[Proof of Theorem~\ref{THM:simplecharacterization}]
We have to show that $S^m({\rm SU}(2))=\Sigma^m({\rm SU}(2))$,
so that theorem would follow from
Theorem~\ref{THM:Gcommutatorcharacterization}.
Both classes $S^m({\rm SU}(2))$ and
$\Sigma^m({\rm SU}(2))$ require
the singular support condition $(y\mapsto R_A(x,y))\subset\{e\}$,
so we do not have to consider this;
moreover, the $x$-dependence of the symbol is not essential here,
and therefore we abbreviate $\sigma_A(l):=\sigma_A(x,l)$.
First, let us show that $\Sigma^m({\rm SU}(2))\subset S^m({\rm SU}(2))$.
Take $\sigma_A\in \Sigma^m({\rm SU}(2))$.
Then also $\sigma_{A_u}\in \Sigma^m({\rm SU}(2))$
(either by the well-known properties of pseudo\-differential 
operators and Theorem \ref{THM:Gcommutatorcharacterization},
or by checking directly
that the definition of the classes $\Sigma_k^m({\rm SU}(2))$
is conjugation-invariant).
Let us define $c_N(B)$ by
\begin{equation*}
  \sigma_{c_N(B)}(l)_{ij} := (i-j)^N\ \sigma_B(l)_{ij}.
\end{equation*}
Now $\sigma_{c_N(A_u)}\in \Sigma^m({\rm SU}(2))$ 
for every $N\in\Bbb Z^+$,
because $\sigma_{A_u}\in \Sigma^m({\rm SU}(2))$ and
$$
  [\sigma_{\partial_0},\sigma_B](l)_{ij} = (i-j)\ \sigma_B(l)_{ij}.
$$
This implies the ``rapid off-diagonal decay'' of $\sigma_{A_u}$:
$$
 \left| 
  \sigma_{A_u}(x,l)_{ij}\right|
  \leq C_{A m N}\ \langle i-j\rangle^{-N}
  \ (1+l)^{m},
$$
implying the norm comparability
\begin{equation}\label{normcomparability}
  \left\| \cdots\sigma_{A_u}(l) \right\|_{op} \sim
  \sup_{i,j}\left| \cdots\sigma_{A_u}(l)_{ij} \right|
\end{equation}
in view of Lemma~\ref{LEM:ell2boundedmatrix} in Appendix.
Moreover,
$\triangle_l^\alpha \partial_x^\beta \sigma_{A_u}
\in \Sigma^{m-|\alpha|}({\rm SU}(2))$ by Corollary
\ref{COR:G-sigma-diff},
so that we obtain the symbol inequalities~\eqref{Sigmainequalities}
from \eqref{norminequality}.
Thereby $\Sigma^m({\rm SU}(2))\subset S^m({\rm SU}(2))$.

Now we have to show that
$S^m({\rm SU}(2))\subset \Sigma^m({\rm SU}(2))$.
Again, we may exploit the norm comparabilities~\eqref{normcomparability}:
thus clearly $S^m({\rm SU}(2))\subset \Sigma_0^m({\rm SU}(2))$.
Consequently, $S^m({\rm SU}(2))\subset \Sigma_k^m({\rm SU}(2))$
for all $k\in\Bbb Z^+$,
due to Lemmas~\ref{LEM:tildeSdiffbehavior} and
\ref{LEM:tildeSproductbehavior}.
\end{proof}

\section{Pseudo-differential operators on manifolds and on $\S3$}
\label{SEC:Pseudos-on-S3}

In this section we discuss how the introduced constructions
are mapped by global diffeomorphisms and give an example
of this in the case of $\SU2$ and $\S3$, proving 
Theorem~\ref{THM:s3-main}.

Let $\Phi:G\to M$ be a diffeomorphism
from a compact Lie group $G$ to a smooth manifold $M$.
Such diffeomorphisms can be obtained for large classes
of compact manifolds by the Poincar\'e conjecture type
results. For example, if $\dim M=3$ it is now known
that such $\Phi$ exists
for any closed simply-connected manifold.

Let us endow $M$ with the natural Lie group structure induced
by $\Phi$,
i.e. with the group multiplication
$((x,y)\mapsto x\cdot y):M\times M\to M$ defined by
$$
  x\cdot y := \Phi\left(\Phi^{-1}(x)\ \Phi^{-1}(y)\right).
$$
Spaces $C^\infty(G)$ and $C^\infty(M)$ are isomorphic
via mappings
\begin{eqnarray*}
  \Phi_\ast: C^\infty(G)\to C^\infty(M), && f\mapsto f_\Phi=f\circ\Phi^{-1},\\
  \Phi^\ast: C^\infty(M)\to C^\infty(G), && g\mapsto g_{\Phi^{-1}}
  = g\circ\Phi.
\end{eqnarray*}
The Haar integral on $M$ is now given by
\begin{equation*}
  \int_M g\ {\rm d}\mu_M \equiv \int_M g\ {\rm d}x
  := \int_G g\circ\Phi\ {\rm d}\mu_G,
\end{equation*}
because for instance
\begin{multline*}
  \int_M g(x\cdot y)\ {\rm d}x
   =  \int_M g(\Phi(\Phi^{-1}(x)\ \Phi^{-1}(y)))\ 
   {\rm d}x = \\
  = \int_G (g\circ\Phi)(\Phi^{-1}(x)\ \Phi^{-1}(y))\ 
   {\rm d}(\Phi^{-1}(x))
   =  \int_G (g\circ\Phi)(z)\ {\rm d}z 
   =  \int_M g(x)\ {\rm d}x.
  \end{multline*}
Moreover,
$\Phi_\ast:C^\infty(G)\to C^\infty(M)$
extends to a linear unitary bijection
$\Phi_\ast:L^2(\mu_G)\to L^2(\mu_M)$:
$$
  \int_M g(x)\ \overline{h(x)}\ {\rm d}x
  = \int_G (g\circ\Phi)\ (\overline{h\circ\Phi})\ {\rm d}\mu_G.
$$
Notice also that there is an isomorphism
\begin{equation*}
  \Phi_\ast: \Rep(G)\to \Rep(M),\quad
  \xi\mapsto \Phi_\ast(\xi)=\xi\circ\Phi
\end{equation*}
of irreducible unitary representations.
Thus $\widehat{G}\cong\widehat{M}$ in this sense.
This immediately implies that the whole construction of
symbols of pseudo-differential operators on $M$ is
equivalent to that on $G$.

Let us now apply this construction to the isomorphism
$\S3\cong{\rm SU}(2)$.
First recall the {quaternion space} ${\mathbb H}$ which is 
the associative $\mathbb R$-algebra
with a vector space basis
$\{{\bf 1},{\bf i}, {\bf j},{\bf k}\}$,
where ${\bf 1}\in\mathbb H$ is the unit and
$$
  {\bf i}^2={\bf j}^2={\bf k}^2
        = -{\bf 1}
        = {\bf ijk}.
$$
Mapping
$
  x=(x_m)_{m=0}^3\mapsto x_0{\bf 1} + x_1{\bf i} + 
  x_2 {\bf j} + x_3 {\bf k}
$
identifies $\mathbb R^4$ with $\mathbb H$.
In particular, the unit sphere 
$\mathbb S^3\subset\mathbb R^4\cong\mathbb H$
is a multiplicative group.
A bijective homomorphism 
$\Phi^{-1}:\S3\to\SU2$ is defined by
\begin{equation*}
  x\mapsto
  \Phi^{-1}(x) = \begin{pmatrix}
        x_0 + {\rm i} x_3 & x_1 + {\rm i} x_2 \\
        -x_1+ {\rm i} x_2 & x_0 - {\rm i} x_3
  \end{pmatrix},
\end{equation*}
and its inverse $\Phi:\SU2\to\S3$ gives rise to
the global quantisation of pseudo-differential operators
on $\S3$ induced by that on
$\SU2$, as shown in the beginning of this section.
This, combined with Theorem~\ref{THM:simplecharacterization},
proves Theorem~\ref{THM:s3-main}.

\section{Appendix on infinite matrices}

In this section we discuss infinite matrices.
The main conclusion that we need is that the 
operator-norm and the $l^\infty$-norm are equivalent for
matrices arising as full symbols of
pseudo-differential operators in $\Psi^m(\SU2)$.

\begin{defn}
Let $\Bbb C^\Bbb Z$ denote the space of complex sequences
$x=(x_j)_{j\in\Bbb Z}$.
A {\it matrix} $A\in\Bbb C^{\Bbb Z\times\Bbb Z}$
is presented as an infinite table
$A = \begin{pmatrix} A_{ij} \end{pmatrix}_{i,j\in\Bbb Z}$.
As usual, we set
$\langle x,y\rangle_{\ell^2}=\sum_{j\in\Bbb Z} x_j\overline{y_j}$,
$\|x\|_{l^2}=\langle x,x\rangle_{\ell^2}^{1/2}$, and
$\|A\|_{l^2\to l^2}=\sup\{\|Ax\|_{l^2}:
\|x\|_{l^2}\leq 1\}$
provided that the sums $(Ax)_i = \sum_{j\in\Bbb Z} A_{ij}\ x_j$
converge absolutely.
For each $k\in\Bbb Z$, let us define
$A(k)\in\Bbb C^{\Bbb Z\times\Bbb Z}$ by
\begin{equation*}
  A(k)_{ij} = \begin{cases} A_{ij}, & {\rm if}\ i-j=k,\\
    0, & {\rm if}\ i-j\not=k. \end{cases}
\end{equation*}
A matrix $A\in\Bbb C^{\Bbb Z\times\Bbb Z}$
will be said to {\it decay (rapidly) off-diagonal} if
\begin{equation}\label{decayoffdiag}
  |A_{ij}| \leq c_{Ar}\ \langle i-j\rangle^{-r}
\end{equation}
for every $i,j\in\Bbb Z$ and $r\in\Bbb N$,
where constants $c_{Ar}<\infty$ depend on $r,A$,
but not on $i,j$.
The set of off-diagonally decaying matrices is 
denoted by $\mathcal D$.
\end{defn}

\begin{lem}\label{LEM:ell2boundedmatrix}
Let $A\in\Bbb C^{\Bbb Z\times\Bbb Z}$
and $\|A\|_{\ell^\infty} = \sup_{i,j\in\Bbb Z} |A_{ij}|$. Then
$$
  \|A\|_{\ell^\infty}\leq \|A\|_{op}.
$$
Moreover, if $|A_{ij}|\leq c\langle i-j\rangle^{-r}$ 
for some $r>1$ then
for $c'=c\sum_{k\in\Bbb Z} \langle k\rangle^{-r}$ we have
$$
  \|A\|_{op}
  \leq c'\ \|A\|_{\ell^\infty}.
$$
\end{lem}

\begin{proof}
Let $\delta_i=(\delta_{ij})_{j\in\Bbb Z}\in\Bbb C^{\Bbb Z}$,
where
$\delta_{ii}=1$ and $\delta_{ij}=0$ if $i\not=j$.
Then
$A_{ij} = \langle A \delta_j, \delta_i\rangle_{\ell^2}$.
The first claim then follows from the Cauchy--Schwarz inequality:
$$
  |A_{ij}|
  =
  \left|\left( A \delta_j, \delta_i\right)_{\ell^2}\right|
  \leq \|A\|_{op}.
$$
Next, since $A=\sum_{k\in\Bbb Z} A(k)$, we get
$$
  \|A\|_{op}
  \leq \sum_{k\in\Bbb Z} \| A(k)\|_{op}
  = \sum_{k\in\Bbb Z} \sup_{j} \left| A(k)_{j+k,j} \right|
  \leq \sum_{k\in\Bbb Z} c \langle k \rangle^{-r}.
$$
From this we directly see that if $\|A\|_{\ell^\infty}\geq 1$
then $\|A\|_{op}\leq c'\ \|A\|_{\ell^\infty}$.
By the linearity of the norms, this concludes the proof.
\end{proof}

\begin{prop}\label{decayingmatrixmultiplication}
Let $A,B\in\mathcal D$. Then $AB\in\mathcal D$.
\end{prop}

\begin{proof}
Matrices $A,B\in\Bbb C^{\Bbb Z\times\Bbb Z}$
in general cannot be multiplied,
but here there is no problem as $A,B\in\mathcal D$, so that
the matrix element $(AB)_{ik}$ is estimated by
$$
  \sum_{j} \left|A_{ij}\right| \left|B_{jk} \right|
  \leq c_{Ar} c_{Bs}
  \sum_{j} \langle i-j\rangle^{-r} \langle j-k\rangle^s
  \stackrel{{\rm Peetre}\; {\rm ineq.}}{\leq}
  c_{Ar} c_{Bs} \sum_{j}
  \langle i-k\rangle^{-r} \langle k-j\rangle^{|r|}
  \langle j-k\rangle^s,
$$
which converges if $|r|+s<-1$.
This shows that $AB\in\mathcal D$.
\end{proof}

Altogether, we obtain the following

\begin{thm}\label{THM:su2-matrices-rapiddecay}
${\mathcal D}\subset{\mathcal L}(\ell^2)$ is a 
unital involutive algebra.
Moreover, for $A\in\mathcal D$,
norms $\|A\|_{op}$ and 
$\|A\|_{\ell^\infty}$
are equivalent.
\end{thm}


\begin{thebibliography}{1}



\bibitem{Agranovich2}
M. S. Agranovich,
{\it Elliptic pseudodifferential operators on a closed curve.}
(Russian) Trudy Moskov. Mat. Obshch.
{\bf 47} (1984), 22--67, 246. 



\bibitem{Beals}
R. Beals,
{\it Characterization of pseudo-differential operators and applications.}
Duke Math. J. {\bf 44} (1977), 45--57.


\bibitem{CoifmanMeyer}
R. R. Coifman and Y. Meyer,
{\it Au-del\`a des op\'erateurs pseudo\-diff\'erentiels.}
Ast\'erisque {\bf 57} (1977), Soci\'et\'e Math. de France.

\bibitem{Co75} H.~O.~Cordes,
  {\it On compactness of commutators of multiplications and convolutions,
    and boundedness of pseudo-differential  operators},
  J. Funct. Anal.
  {\bf 18} (1975), 115--131.


\bibitem{Dunau}
J. Dunau,
{\it Fonctions d'un operateur elliptique sur une variete compacte.}
J. Math. Pures et Appl. {\bf 56} (1977), 367--391.


\bibitem{GK97}
F. Geshwind and N. H. Katz, 
{\it Pseudodifferential operators on ${\rm SU}(2)$},
J. Fourier Anal. Appl. {\bf 3} (1997), 193--205. 

\bibitem{Gl04}
P. Glowacki, 
{\it A symbolic calculus and $L\sp 2$-boundedness on 
nilpotent Lie groups},
J. Funct. Anal. {\bf 206} (2004), 233--251. 


 \bibitem{HR} E. Hewitt and K. A. Ross,
       {\em Abstract harmonic analysis I, II.}
       Springer-Verlag, 1963, 1970.


\bibitem{Hormander}
L. H\"ormander,
{\it Pseudo-differential operators.}
Comm. Pure Appl. Math. {\bf 18} (1965), 501--517.

\bibitem{Hormander2}
L. H\"ormander,
{\it The Analysis of Linear Partial Differential Operators III.}
Berlin: Springer-Verlag, 1985.

\bibitem{How84}
R. Howe, 
{\it A symbolic calculus for nilpotent groups}, 
Operator algebras and group representations, 
Vol. I (Neptun, 1980), 254--277, 
Monogr. Stud. Math., 17, Pitman, Boston, MA, 1984. 


\bibitem{KohnNirenberg}
J. J. Kohn and L. Nirenberg,
{\it On the algebra of pseudo-differential operators.}
Comm. Pure Appl. Math. {\bf 18} (1965), 269--305.

\bibitem{McLean}
W. McLean,
{\it Local and global descriptions of periodic 
pseudo\-differential operators.}
Math. Nachr. {\bf 150} (1991), 151--161.

\bibitem{MeSh87}
G. A. Meladze and M. A. Shubin, 
{\it A functional calculus of pseudodifferential 
operators on unimodular Lie groups}. 
J. Soviet Math. {\bf 47} (1989), 2607--2638.

\bibitem{Me83}
A. Melin, 
{\em Parametrix constructions for right invariant 
differential operators on nilpotent groups},
Ann. Global Anal. Geom. {\bf 1} (1983), 79--130. 





\bibitem{RT07}
M. Ruzhansky and V. Turunen, 
{\it On the Fourier analysis of operators on the torus}, 
Modern trends in pseudo-differential operators, 87-105, 
Oper. Theory Adv. Appl., 172, Birkh\"auser, Basel, 2007.

\bibitem{RT08}
M. Ruzhansky and V. Turunen, 
{\it Quantization of pseudo-differential operators on
 the torus}, arXiv:0805.2892v1

\bibitem{RT-book}
M. Ruzhansky and V. Turunen, 
{\it Pseudo-differential operators and symmetries},
monograph in preparation, to appear in Birkh\"auser.

\bibitem{Sa97} Yu. Safarov,
\textit{Pseudodifferential operators and 
linear connections.} Proc. London Math. Soc.
\textbf{74} (1997), 379--416.


\bibitem{SV02}
J. Saranen and G. Vainikko,
{\it Periodic integral and pseudodifferential 
equations with numerical approximation.} 
Springer Monographs in Mathematics.
Springer-Verlag, Berlin, 2002.

\bibitem{SW87}
J. Saranen and W. L. Wendland,
{\it The Fourier series representation of 
pseudodifferential operators on closed curves.} 
Complex Variables Theory Appl. {\bf 8} (1987), 55--64.

\bibitem{Segal1}
I. E. Segal, 
{\it An extension of Plancherel's formula to 
separable unimodular groups.}
Ann. Math. {\bf 52} (1950), 272--292.


\bibitem{Sh05} V. A. Sharafutdinov,
\textit{Geometric symbol calculus for 
pseudodifferential operators. I} 
[Translation of Mat. Tr. \textbf{7} (2004), 159--206].
Siberian Adv. Math. \textbf{15} (2005), 81--125.

\bibitem{She75}
T. Sherman, {\it Fourier analysis on the sphere}, 
Trans. Amer. Math. Soc. {\bf 209} (1975), 1--31.

\bibitem{She90}
T. Sherman, 
{\it The Helgason Fourier transform for compact 
Riemannian symmetric spaces of rank one},
Acta Math. {\bf 164} (1990), 73--144. 



\bibitem{St72}
R. Strichartz,
{\it Invariant pseudo-differential operators on a Lie group.}
Ann. Scuola Norm. Sup. Pisa {\bf 26} (1972), 587--611.

\bibitem{St88}
R. Strichartz, 
{\it Local harmonic analysis on spheres.}
J. Funct. Anal. {\bf 77} (1988), 403--433. 




\bibitem{Taylor84}
M. E. Taylor,
{\it Noncommutative microlocal analysis.}
Mem. Amer. Math. Soc. {\bf 52} (1984), No. 313.

\bibitem{Taylor86}
M. E. Taylor,
{\it Noncommutative harmonic analysis.}
Mathematical Surveys and Monographs, Vol. 22,
Amer. Math. Soc., 1986.

\bibitem{Taylor97}
M. E. Taylor,
{\it Beals--Cordes -type characterizations of pseudodifferential
operators.}
Proc. Amer. Math. Soc. {\bf 125} (1997), 6, 1711--1716.


\bibitem{Turunen00}
V. Turunen,
{\it Commutator characterization of periodic 
pseudo\-differential operators.}
Z. Anal. Anw. {\bf 19} (2000), 95--108. 

\bibitem{Turunen01}
V. Turunen,
{\it Pseudodifferential calculus on compact Lie groups.},
Helsinki Univ. Techn. Inst. Math. Research Report A431. 2001.


\bibitem{TurunenVainikko}
V. Turunen and G. Vainikko,
{\it On symbol analysis of periodic pseudo\-differential operators.}
Z. Anal. Anw. {\bf 17} (1998), 9--22. 


\bibitem{Vilenkin}
N. Vilenkin,
{\it Special functions and the theory of group representations.}
Trans. Math. Monographs, Vol. 22, Amer. Math. Soc., 1968.



\bibitem{Wi80}
H. Widom,
\textit{A complete symbolic calculus for 
pseudodifferential operators}.
Bull. Sci. Math. \textbf{104} (1980), 19--63.



\end{thebibliography}
\end{document}